\documentclass[preprint,11pt]{elsarticle}
\usepackage[T1]{fontenc}
\usepackage[letterpaper,margin=1.14in]{geometry}
\usepackage{amsmath,amssymb,amsthm,microtype,mathrsfs}
\usepackage[hidelinks]{hyperref}
\hypersetup{pdftitle={Schauder bases and local unconditional structure in Hardy spaces of Dirichlet series}}
\patchcmd{\pprintMaketitle}{\hrule\vskip12pt}{\vskip12pt}{}{}
\patchcmd{\pprintMaketitle}{\hrule\vskip12pt}{\vskip12pt}{}{}
\patchcmd{\pprintMaketitle}{\footnotesize\itshape\elsaddress\par\vskip36pt}{}{}{}
\allowdisplaybreaks[2]
\newtheorem{theorem}{Theorem}[section]
\newtheorem{lemma}[theorem]{Lemma}
\newtheorem{proposition}[theorem]{Proposition}
\theoremstyle{definition}
\newtheorem{definition}[theorem]{Definition}
\theoremstyle{plain}
\numberwithin{equation}{section}
\DeclareMathOperator{\lust}{lust}
\begin{document}
\begin{frontmatter}
\title{Schauder bases and local unconditional structure in Hardy spaces of Dirichlet series}
\author{Tengfei Ma}
\author{Yufeng Lu}
\author{Chao Zu\corref{corresponding}}
\cortext[corresponding]{Corresponding author.}
\begin{abstract}
We prove that among the Hardy spaces \(\mathcal H^p\) of Dirichlet
series with \(1\le p<\infty\), only \(\mathcal H^2\) has local unconditional
structure. Consequently, \(\mathcal H^p\) admits an unconditional basis
exactly when \(p=2\). At the endpoint \(p=1\), we also construct a Schauder
basis of Dirichlet polynomials and give an explicit formula for each basis
vector. Together, these results answer the basis questions posed by Aleman,
Olsen and Saksman, as well as the question of Saksman and Seip concerning
an unconditional basis for \(\mathcal H^1\).
\end{abstract}
\end{frontmatter}
\pagestyle{plain}\thispagestyle{plain}

\section{Introduction}

Dirichlet series \(f(s)=\sum_{n=1}^{\infty}a_n n^{-s}\) have a long history in analytic number theory \cite[Ch.~I, \S1]{HardyRiesz}. Their modern function-theoretic study is closely connected with a fundamental observation of Bohr: the multiplicative structure of the positive integers allows a Dirichlet series to be viewed as a power series in infinitely many variables. If \(p_j\) is the \(j\)-th prime and
\[
 n=\prod_{j\ge1}p_j^{\alpha_j(n)},
\]
then the substitution \(z_j=p_j^{-s}\) associates \(n^{-s}\) with the monomial
\[
 z^{\alpha(n)}=\prod_{j\ge1}z_j^{\alpha_j(n)}.
\]
Bohr's correspondence between formal Dirichlet series and formal power series is
\[
 f(s)=\sum_{n=1}^{\infty}a_n n^{-s}
 \quad\longleftrightarrow\quad
 \sum_{n=1}^{\infty}a_n z^{\alpha(n)} .
\]
This correspondence is more than a formal change of variables: it identifies the analysis of Dirichlet series with Fourier analysis of functions on the infinite-dimensional torus. It is the basic link between the arithmetic structure of the integers and function theory in infinitely many variables; see \cite[Chapters~3 and~11]{DGMS} and \cite{Survey}. Bohr's viewpoint also underlies the classical work of Bohnenblust and Hille on absolute convergence, and provides an early connection between Dirichlet series and holomorphic functions on polydiscs \cite{BohnenblustHille}.

Although the underlying correspondence is classical, the systematic study of Hardy spaces of Dirichlet series from the viewpoint of modern functional analysis began with the work of Hedenmalm, Lindqvist, and Seip \cite{HLS}. Their 1997 paper was motivated in part by completeness and Riesz basis questions for systems of dilated functions \(\{\varphi(nx)\}_{n\ge1}\) in \(L^2(0,1)\). They introduced the Hilbert space
\[
 \mathcal H^2=\left\{\sum_{n\ge1}a_n n^{-s}:\ \sum_{n\ge1}|a_n|^2<\infty\right\},
\]
whose elements are analytic in the half-plane \(\operatorname{Re}s>1/2\), and used Bohr's correspondence to identify it with the Hardy space on \(\mathbb T^\infty\). This identification brought the theory of Dirichlet series into direct contact with harmonic analysis on compact groups and the function theory of infinitely many complex variables.

Bayart \cite{Bayart} subsequently developed the corresponding \(\mathcal H^p\) spaces for \(p\ge 1\) and studied, among other things, composition operators and a Poisson-type semigroup. Under the Bohr correspondence, the space \(\mathcal H^p\) is identified isometrically with the analytic Hardy space \(H^p(\mathbb T^\infty)\), the closed subspace of \(L^p(\mathbb T^\infty)\) generated by the analytic monomials. This model has made it possible to bring tools from several complex variables, harmonic analysis, and probability into the study of Dirichlet series. The resulting theory includes work on Fourier multipliers \cite{AAOS}, convergence and boundary behaviour \cite{HS}, composition operators \cite{Bayart}, and a range of embedding, interpolation, and operator-theoretic questions; see also the monograph \cite{DGMS} and the survey of open problems \cite{Survey}. The interaction with analytic number theory is equally important, since the multiplicative structure of the coefficients is reflected in the geometry and norm estimates of these spaces.

The natural ordering of the Dirichlet monomials is the first candidate for a basis. In \(\mathcal H^2\), these monomials form an orthonormal basis. More generally, Helson's theorem \cite{Helson} implies that the partial sum operators
\[
 S_N\!\left(\sum_{n\ge1}a_n n^{-s}\right)=\sum_{n\le N}a_n n^{-s}
\]
are uniformly bounded on \(\mathcal H^p\) for \(1<p<\infty\); consequently, the Dirichlet monomials form a Schauder basis in this range. Their unconditional behaviour is different: when \(p\ne2\), the canonical basis is not unconditional \cite[Section~4]{CLM}. This leaves open the more substantial question of whether some other choice of basis could be unconditional.

At the endpoint \(p=1\), even the natural ordering ceases to give a Schauder basis. Indeed, the upper estimate of Bondarenko, Brevig, Saksman, and Seip \cite{BBSS} and the lower estimate of Andersson \cite{Andersson} yield
\[
 \|S_N\|_{\mathcal H^1\to\mathcal H^1}\asymp\frac{\log N}{\log\log N},
 \qquad N\longrightarrow\infty ,
\]
where \(\asymp\) means that each side bounds the other up to an absolute constant.
In particular, the partial sum operators are not uniformly bounded. This rules out the canonical ordering, but it does not rule out a Schauder basis formed from different Dirichlet polynomials. The existence of such a basis is therefore a genuinely separate question.

These basis questions were raised by Aleman, Olsen, and Saksman \cite[Questions~1 and~2]{AAOS}, who asked whether \(\mathcal H^1\) admits a Schauder basis and whether \(\mathcal H^p\) admits an unconditional basis for \(1<p<\infty\), \(p\ne2\). Saksman and Seip \cite[Problems~2.5 and~2.6]{Survey} also asked whether \(\mathcal H^1\) has an unconditional basis. We approach these questions through the notion of local unconditional structure introduced by Gordon and Lewis \cite{GL}. Our first result gives a complete classification.

\begin{theorem}\label{thm:unc}
For \(1\le p<\infty\), the Hardy space \(\mathcal H^p\) has local unconditional structure if and only if \(p=2\). In particular, \(\mathcal H^p\) admits an unconditional basis if and only if \(p=2\).
\end{theorem}

Thus, outside the Hilbert space case, no choice of basis can make the expansions unconditional. At the same time, the failure of local unconditional structure does not preclude the existence of a Schauder basis. We prove that this is precisely what happens at the endpoint.

\begin{theorem}\label{thm:basis}
The space \(\mathcal H^1\) admits an explicit Schauder basis of Dirichlet polynomials. Every Dirichlet polynomial has a finite expansion in this basis.
\end{theorem}

The basis vectors and their coefficient functionals are given explicitly in formulas \eqref{eq:explicit-vector} and \eqref{eq:explicit-dual}, respectively. In particular, the theorem gives a concrete basis in the original prime coordinates of the Dirichlet series, rather than only an abstract existence result.

Section~\ref{sec:pre} introduces the notation and preliminary facts on bases and local unconditional structure. Both proofs use the finite-dimensional constructions developed in Section~\ref{sec:finite}. We prove the classification in Section~\ref{sec:unc-proof} and construct the Schauder basis of \(\mathcal H^1\) in Section~\ref{sec:basis-proof}, where we give its explicit representation in the original Dirichlet series variables.

\section{Preliminaries}\label{sec:pre}

\subsection{Hardy spaces of Dirichlet series}

Throughout, we assume \(1\leq p<\infty\), unless a different range is explicitly stated. All spaces are over the complex field unless otherwise specified.
We write \(u\lesssim v\), or equivalently \(v\gtrsim u\), when \(u\le Cv\)
for an absolute constant \(C\). Subscripts on comparison signs indicate the
parameters on which the constants may depend.

\smallskip\noindent\textbf{Hardy spaces on the infinite torus.}
Let \(\mathbb T=\{z\in\mathbb C:|z|=1\}\) and \(\mathbb T^\infty=\prod_{j\geq1}\mathbb T\).
We write \(m\) for normalized Haar measure on each torus. Denote by \(\mathbb Z^{(\mathbb N)}\) the set of finitely supported integer sequences and by \(\mathbb N_0^{(\mathbb N)}\) its subset consisting of sequences with nonnegative entries. For \(\nu\in\mathbb Z^{(\mathbb N)}\), write
\[
z^\nu=\prod_{j\geq1}z_j^{\nu_j}.
\]
The Fourier coefficient of \(f\in L^p(\mathbb T^\infty,m)\) at \(\nu\) is
\[
\widehat f(\nu)
=\int_{\mathbb T^\infty}f(z)z^{-\nu}\,dm(z).
\]

The analytic Hardy space on the infinite torus is defined by
\[
H^p(\mathbb T^\infty)
=
\left\{
f\in L^p(\mathbb T^\infty,m):
\widehat f(\nu)=0
\text{ whenever }\nu\notin\mathbb N_0^{(\mathbb N)}
\right\},
\]
with the inherited \(L^p\)-norm. Equivalently, \(H^p(\mathbb T^\infty)\) is the closed linear span in \(L^p(\mathbb T^\infty)\) of the analytic monomials \(z^\alpha\), where \(\alpha\in\mathbb N_0^{(\mathbb N)}\); see \cite[Chapter~5, (5.5), and Theorem~11.10]{DGMS}. We use the notation
\[
|\alpha|=\sum_{j\geq1}\alpha_j.
\]
The corresponding spaces on \(\mathbb T^d\) are defined in the same way, with Fourier indices in \(\mathbb Z^d\) and analytic indices in \(\mathbb N_0^d\).

\smallskip\noindent\textbf{Bohr transform.}
Let \(p_1<p_2<\cdots\) denote the prime numbers. Every positive integer \(n\) has a unique factorization
\[
n=\prod_{j\geq1}p_j^{\alpha_j(n)},
\]
where \(\alpha(n)\in\mathbb N_0^{(\mathbb N)}\), and we set \(\alpha(1)=0\). The Bohr transform is the linear correspondence between formal power series in countably many variables and formal Dirichlet series given by
\begin{equation}\label{eq:bohr-transform}
\mathcal B\left(\sum_{\alpha\in\mathbb N_0^{(\mathbb N)}}c_\alpha z^\alpha\right)
=
\sum_{n\geq1}c_{\alpha(n)}n^{-s}.
\end{equation}
We follow the convention of \cite[Chapter~3, (3.2)]{DGMS}.

\smallskip\noindent\textbf{Hardy spaces of Dirichlet series.}
We define \(\mathcal H^p\) as the image of \(H^p(\mathbb T^\infty)\) under the Bohr transform, equipped with the transported norm:
\[
\mathcal H^p=\{\mathcal Bf:f\in H^p(\mathbb T^\infty)\},
\qquad
\|\mathcal Bf\|_{\mathcal H^p}
=
\|f\|_{L^p(\mathbb T^\infty,m)}.
\]
Consequently, the Bohr transform is an isometric isomorphism
\[
\mathcal B:H^p(\mathbb T^\infty)\longrightarrow \mathcal H^p;
\]
see \cite[Chapter~11, (11.6)]{DGMS}. In particular, if
\[
\mathcal Bf(s)=\sum_{n\geq1}a_n n^{-s},
\]
then
\[
a_n=\widehat f(\alpha(n)).
\]

\begin{proposition}[Dirichlet polynomials]\label{prop:dirichlet-means}
Dirichlet polynomials are dense in \(\mathcal H^p\). Moreover, for every Dirichlet polynomial
\[
P(s)=\sum_{n=1}^N a_n n^{-s},
\]
we have
\[
\|P\|_{\mathcal H^p}
=
\lim_{R\to\infty}
\left(
\frac{1}{2R}\int_{-R}^R|P(it)|^p\,dt
\right)^{1/p}.
\]
\end{proposition}

The density assertion and the norm identity follow from \cite[Theorems~11.9 and 11.10]{DGMS}.

\subsection{Coordinate truncations on the infinite torus}\label{sec:coordinate-sections}

For a formal power series \(F\) in \(z_1,z_2,\ldots\) and an integer \(m\geq1\), define
\begin{equation}\label{eq:coordinate-section}
(A_mF)(z_1,\ldots,z_m)=F(z_1,\ldots,z_m,0,0,\ldots).
\end{equation}
For \(f\in H^p(\mathbb T^\infty)\), we apply \eqref{eq:coordinate-section} to its formal Fourier series.
We identify a function on \(\mathbb T^m\) with its extension to \(\mathbb T^\infty\) that depends only on the first \(m\) coordinates.

\begin{lemma}[{\cite[Theorem~2.1]{BBSS}}]\label{lem:prime-sections}
Let \(1\leq p<\infty\), \(m\geq1\), and \(f\in H^p(\mathbb T^\infty)\). Then \(A_mf\in H^p(\mathbb T^m)\) and
\[
\|A_mf\|_{H^p(\mathbb T^m)}\leq\|f\|_{H^p(\mathbb T^\infty)}.
\]
\end{lemma}

\subsection{Bases and local unconditional structure}\label{sec:bases-lust}

\begin{definition}\label{def:bases}
A sequence \((x_j)_{j\geq1}\) in a Banach space \(X\) is called a Schauder basis if for every \(x\in X\) there are unique coefficients \(a_1,a_2,\ldots\) such that
\[
\lim_{N\to\infty}
\left\|x-\sum_{j=1}^N a_jx_j\right\|=0.
\]
The basis constant is \(\sup_{N\geq1}\|P_N\|\), where \(P_Nx=\sum_{j=1}^N a_jx_j\) is the \(N\)th partial sum of the expansion of \(x\).

The basis is unconditional if, for every \(x\in X\), every rearrangement of its basis expansion converges in norm to \(x\).
Equivalently, there exists \(K\geq1\) such that, for every \(N\geq1\), all coefficients \(a_1,\ldots,a_N\), and all \(\lambda_1,\ldots,\lambda_N\) with \(|\lambda_j|=1\),
\[
\left\|\sum_{j=1}^N\lambda_ja_jx_j\right\|
\leq
K\left\|\sum_{j=1}^Na_jx_j\right\|.
\]
The unconditional basis constant is the infimum of all admissible \(K\). If this infimum equals \(1\), the basis is called \(1\)-unconditional.
\end{definition}

We use the following notion of local unconditional structure.

\begin{definition}\label{def:lust}
A Banach space \(X\) has local unconditional structure if there exists a constant \(C\geq1\) such that, for every finite-dimensional subspace \(F\subset X\), there are a Banach space \(W\) with a \(1\)-unconditional basis and bounded linear operators
\[
u:F\longrightarrow W,
\qquad
v:W\longrightarrow X
\]
satisfying
\[
v(u(x))=x\quad(x\in F),
\qquad
\|u\|\,\|v\|\leq C.
\]
We denote by \(\lust(X)\) the infimum of all admissible constants \(C\), with the convention that \(\lust(X)=\infty\) if no such constant exists.
\end{definition}

This is the Gordon--Lewis notion of local unconditional structure \cite[\S3]{GL}, with the normalization used in \cite[Definition~2.1]{KT}. It is also referred to as GL-local unconditional structure in \cite[Appendix]{Johnson}.

Local unconditional structure is inherited by complemented subspaces. The estimates below follow from the remark after Lemma~3.1 and from Lemma~3.2 in \cite{GL}.

\begin{lemma}\label{lem:lust-transfer}
Let \(X\) be a Banach space.
\begin{enumerate}
\item[\textup{(i)}] If \((x_j)\) is an unconditional basis of \(X\), then \(\lust(X)\leq K\), where \(K\) is the unconditional basis constant of \((x_j)\).
\item[\textup{(ii)}] Suppose the identity on a nonzero Banach space \(Y\) factors as \(I_Y=ri\), with bounded linear maps \(i:Y\to X\) and \(r:X\to Y\). Then
\begin{equation}\label{eq:transfer}
\lust(Y)\leq\|i\|\,\|r\|\,\lust(X).
\end{equation}
When \(Y\) is a complemented subspace of \(X\) with projection \(P:X\to Y\), the inclusion \(i:Y\hookrightarrow X\) and \(r=P\) give \(\lust(Y)\leq\|P\|\,\lust(X)\).
\end{enumerate}
\end{lemma}

A sequence \((D_j)_{j\geq1}\) of finite-dimensional subspaces of \(X\) is called a finite-dimensional decomposition (FDD) if every \(x\in X\) has a unique norm-convergent representation
\[
x=\sum_{j\geq1}d_j,
\qquad d_j\in D_j.
\]
The FDD is unconditional if there exists \(K\geq1\) such that, for every \(N\), every choice of \(d_j\in D_j\), and every choice of signs \(\varepsilon_j\in\{-1,1\}\),
\[
\left\|\sum_{j=1}^N\varepsilon_jd_j\right\|
\leq
K\left\|\sum_{j=1}^Nd_j\right\|.
\]
We use the sign convention for unconditional decompositions adopted in \cite[\S2]{KT}.

The \(\ell^1\)-direct sum of a sequence of Banach spaces \((V_j)_{j\geq1}\) is
\[
\left(\bigoplus_{j\geq1}V_j\right)_1
=
\left\{
(v_j):v_j\in V_j,\ 
\sum_{j\geq1}\|v_j\|_{V_j}<\infty
\right\},
\]
equipped with the norm
\[
\|(v_j)\|=\sum_{j\geq1}\|v_j\|_{V_j}.
\]
When all the summands coincide with \(X\), we write \(\ell^1(X)\).

A Banach space \(X\) has cotype \(r\), where \(2\leq r<\infty\), if there is a constant \(C\) such that, for every finite family \(x_1,\ldots,x_N\in X\),
\begin{equation}\label{eq:cotype}
\left(\sum_{j=1}^N\|x_j\|^r\right)^{1/r}
\leq
C\left(\mathbb E\left\|\sum_{j=1}^N\varepsilon_jx_j\right\|^2\right)^{1/2},
\end{equation}
where \((\varepsilon_j)\) are independent Rademacher random variables. For \(1\leq p<\infty\), the space \(L^p(\mu)\) has cotype \(\max(2,p)\), with a constant depending only on \(p\); see \cite[\S1.e, pp.~72--73]{LT}.

\subsection{Fourier analysis on finite groups and probabilistic estimates}

Let \(G\) be a finite abelian group, equipped with normalized counting measure. Its dual group, denoted by \(\widehat G\), consists of all homomorphisms \(\chi:G\to\mathbb T\), with pointwise multiplication. For \(f:G\to\mathbb C\), define
\begin{equation}\label{eq:finite-fourier}
\widehat f(\chi)
=
\frac1{|G|}\sum_{a\in G}f(a)\overline{\chi(a)},
\qquad
L_A^p(G)
=
\operatorname{span}\{\chi:\chi\in A\},
\quad A\subset\widehat G.
\end{equation}
The space \(L_A^p(G)\) is equipped with the norm inherited from \(L^p(G)\).

The characters form an orthonormal basis of \(L^2(G)\). In particular, the Fourier inversion formula gives
\[
f=\sum_{\chi\in\widehat G}\widehat f(\chi)\chi;
\]
see \cite[Chapter~7, Theorem~2.7]{SteinShakarchi}. Since every point of \(G\) has positive measure, elements of \(L^p(G)\) have well-defined point values. We write \(1_A\) for the indicator function of a set \(A\).

\begin{lemma}[{\cite[Lemma~2.2]{Dirksen}}]\label{lem:unit-moment}
Let \(v_1,\ldots,v_N\) belong to a complex Hilbert space, and let \(\eta_1,\ldots,\eta_N\) be independent complex random variables satisfying
\[
\mathbb E\eta_j=0,
\qquad
|\eta_j|=1\quad\mbox{almost surely}.
\]
Then
\begin{equation}\label{eq:unit-moment}
\mathbb E\left\|\sum_{j=1}^N\eta_jv_j\right\|
\asymp
\left(\sum_{j=1}^N\|v_j\|^2\right)^{1/2},
\end{equation}
with absolute constants.
\end{lemma}

\begin{lemma}\label{lem:martingale-bounds}
Let \(M=(M_j)_{j=0}^b\) be a complex-valued martingale, and put
\[
M^*=\max_{0\leq j\leq b}|M_j|.
\]
If \(M_b\in L^2\), then
\begin{equation}\label{eq:martingale-max}
\|M^*\|_2\leq2\|M_b\|_2.
\end{equation}
Moreover, if \(M\) is an integrable martingale with \(M_0=0\), then
\begin{equation}\label{eq:martingale-square}
\mathbb E M^*
\asymp
\mathbb E\left(\sum_{j=1}^b|M_j-M_{j-1}|^2\right)^{1/2}.
\end{equation}
\end{lemma}

\begin{proof}
By conditional Jensen's inequality, \((|M_j|)_{j=0}^b\) is a nonnegative submartingale. Applying Doob's maximal inequality \cite[Thm.~3.4, p.~317]{Doob} to this process with \(\alpha=2\) gives
\[
\mathbb E (M^*)^2\leq 4\,\mathbb E|M_b|^2.
\]

Now assume \(M_0=0\) and write \(M=U+iV\). For a real martingale \(X=U\) or \(V\), extended constantly after time \(b\), the Burkholder--Davis--Gundy inequality \cite[Thm.~1.1]{BDG} with \(\Phi(t)=t\) gives
\[
\mathbb E X^*\asymp\mathbb E S(X),\qquad
S(X)=\left(\sum_{j=1}^b|X_j-X_{j-1}|^2\right)^{1/2}.
\]
Since \(M^*\leq U^*+V^*\leq 2M^*\), we obtain
\[
\mathbb E M^*\asymp\mathbb E\bigl(S(U)+S(V)\bigr)
\asymp\mathbb E\left(S(U)^2+S(V)^2\right)^{1/2}
=\mathbb E\left(\sum_{j=1}^b|M_j-M_{j-1}|^2\right)^{1/2}.
\qedhere
\]
\end{proof}

\section{Finite-dimensional embeddings and basis estimates}\label{sec:finite}

\subsection{Complemented embeddings of \texorpdfstring{\(L_A^p(G)\)}{L\string^p\string_A(G)} into \texorpdfstring{\(H^p(\mathbb T^\infty)\)}{H\string^p(T\string^infinity)}}

The following lemma uses Jeroslow's observation that each vertex of
\(\{0,1\}^b\) can be excluded by a single affine inequality satisfied by
all the other vertices \cite[p.~119, following (1)]{Jeroslow}.

\begin{lemma}[Affine inequalities on the cube]\label{lem:cube-affine}
Let \(b\) be a positive integer and let \(D\subset\{0,1\}^b\).
There exist an integer \(m\ge1\), a matrix \(B\in\mathbb Z^{m\times b}\),
and a vector \(c\in\mathbb Z^m\) such that, for \(\Phi(\alpha)=c+B\alpha\),
\begin{equation}\label{eq:cube-affine-criterion}
 D=\{\alpha\in\mathbb Z^b:\ \alpha\ge0,\ \Phi(\alpha)\ge0\},
\end{equation}
where both inequalities are coordinatewise.
\end{lemma}

\begin{proof}
Let \(V=\{0,1\}^b\setminus D\) and \(m=b+|V|\). Define
\(B\in\mathbb Z^{m\times b}\) and \(c\in\mathbb Z^m\) by
\[
 B=\begin{pmatrix}
 -I_b\\[2pt]
 (1-2v_j)_{v\in V,\,1\le j\le b}
 \end{pmatrix},
 \qquad
 c=\begin{pmatrix}
 (1,\ldots,1)^{\mathrm T}\\[2pt]
 (\sum_{j=1}^b v_j-1)_{v\in V}
 \end{pmatrix}.
\]
Use the same ordering of
\(V\) for the last \(|V|\) rows of \(B\) and \(c\).

The first \(b\) entries of \(c+B\alpha\) are \(1-\alpha_i\),
\(1\le i\le b\). For \(\alpha\in\mathbb Z^b\), these entries are
nonnegative and \(\alpha\ge0\) precisely when \(\alpha\in\{0,1\}^b\).
At a vertex \(\alpha\in\{0,1\}^b\), the row indexed by \(v\in V\) gives
\[
 \begin{aligned}
 (c+B\alpha)_v
 &=\sum_{j=1}^b v_j-1+\sum_{j=1}^b(1-2v_j)\alpha_j\\
 &=\sum_{j:\,v_j=0}\alpha_j+\sum_{j:\,v_j=1}(1-\alpha_j)-1\\
 &=|\{j:\alpha_j\ne v_j\}|-1.
 \end{aligned}
\]
This entry is negative precisely when \(\alpha=v\). Hence, for every \(\alpha\in\mathbb Z^b\),
\[
 \alpha\ge0 \ \text{and}\ c+B\alpha\ge0
 \quad\Longleftrightarrow\quad \alpha\in\{0,1\}^b\setminus V
 \quad\Longleftrightarrow\quad \alpha\in D.
 \qedhere
\]
\end{proof}

The next theorem gives an almost isometric embedding of \(L_A^p(G)\) into
\(H^p(\mathbb T^\infty)\) with complemented range.
The proof draws on the isometric embedding of trigonometric polynomials
used by Konyagin, Queff\'elec, Saksman and Seip
\cite[Lemma~2.5 and the proof of Lemma~2.6]{KQSS} and on the Haar
averaging formula for Fourier projections used by Brevig, Ortega-Cerd\`a and Seip
\cite[\S2.1, (2.2)]{BOS}.

\begin{theorem}\label{thm:retract}
Let \(1\le p<\infty\), let \(G\) be a finite abelian group, and fix \(A\subset\widehat G\) and \(\varepsilon>0\). There are bounded linear operators
\[
 J:L^p(G)\longrightarrow L^p(\mathbb T^\infty),\qquad
 R:L^p(\mathbb T^\infty)\longrightarrow L^p(G)
\]
such that
\begin{equation}\label{eq:retraction-bounds}
 RJ=I,\qquad \|J\|\,\|R\|\le 1+\varepsilon,
\end{equation}
and \(J(L^p_A(G))\subset H^p(\mathbb T^\infty)\),
\(R(H^p(\mathbb T^\infty))\subset L^p_A(G)\).
\end{theorem}

\begin{proof}
\textbf{The trivial group.}
If \(|G|=1\), then \(L^p(G)=\mathbb C\) and \(\widehat G=\{1\}\).
Take \(u=1\) if \(A=\{1\}\) and \(u(z)=\overline{z_1}\) if \(A=\varnothing\). Set
\[
 Jc=cu,\qquad RF=\int_{\mathbb T^\infty}F\,\overline u\,dm .
\]
Since \(|u|=1\), we have \(RJ=I_{\mathbb C}\) and
\(\|J\|=\|R\|=1\).
If \(A=\{1\}\), the range of \(J\) consists of constant functions.
If \(A=\varnothing\), then \(L_A^p(G)=\{0\}\) and
\(RF=\widehat F(-1,0,\ldots)=0\) for every \(F\in H^p(\mathbb T^\infty)\).

Assume henceforth that \(|G|>1\).

\textbf{Step 1. Construct an almost isometric embedding of \(L^p(G)\) by polynomials.}
By the structure theorem for finite abelian groups
\cite[Theorem~13.4]{Judson}, we may identify \(G\) with \(\prod_{k=1}^d\mathbb Z/N_k\mathbb Z\),
where \(N_k\ge2\) is the order of the \(k\)th cyclic factor.
Write \(a=(a_1,\ldots,a_d)\in G\), with \(a_k\in\mathbb Z/N_k\mathbb Z\),
and put \(\omega_k=e^{2\pi i/N_k}\). For \(m\in\mathbb Z\), let
\([m]_{N_k}=m+N_k\mathbb Z\) denote its residue class modulo \(N_k\). We set
\(\omega_k^{[m]_{N_k}}=\omega_k^m\), which is well defined since \(\omega_k^{N_k}=1\).
For each positive integer \(r\), write \(w=(w_{k,j})\in\mathbb T^{dr}\),
with coordinates indexed by \((k,j)\), \(1\le k\le d\), \(1\le j\le r\), and set
\begin{equation}\label{eq:character-kernel}
 H_a^{(r)}(w)=\prod_{k=1}^d\prod_{j=1}^r
 \frac{1+\omega_k^{-a_k}w_{k,j}}{\|1+z\|_{L^p(\mathbb T)}}.
\end{equation}
Define \(S:L^p(G)\to L^p(\mathbb T^{dr})\) by
\begin{equation}\label{eq:finite-embedding}
 (Sf)(w)=|G|^{-1/p}\sum_{a\in G}f(a)H_a^{(r)}(w).
\end{equation}
By rotation invariance and the product formula for Haar measure,
\[
 \|H_a^{(r)}\|_{L^p(\mathbb T^{dr})}^p
 =\prod_{k,j}\frac{\|1+\omega_k^{-a_k}z\|_{L^p(\mathbb T)}^p}
                     {\|1+z\|_{L^p(\mathbb T)}^p}=1 .
\]
The continuous function \(|H_a^{(1)}|\) has its unique maximum at
\((\omega_k^{a_k})_{k=1}^d\), so for \(a\ne b\),
\[
 \int_{\mathbb T^d}|H_a^{(1)}H_b^{(1)}|^{p/2}\,dm
 =1-\tfrac12\big\|\,|H_a^{(1)}|^{p/2}-|H_b^{(1)}|^{p/2}\,\big\|_{L^2(\mathbb T^d)}^2<1 .
\]
Since \(G\) is finite, we may set
\[
 \kappa=\max_{\substack{a,b\in G\\a\ne b}}
 \int_{\mathbb T^d}|H_a^{(1)}H_b^{(1)}|^{p/2}\,dm<1 .
\]

The defining product factors as
\[
 H_a^{(r)}(w)=\prod_{j=1}^r H_a^{(1)}(w_{1,j},\ldots,w_{d,j}).
\]
Integrating separately in the \(r\) groups of variables gives, for \(a\ne b\),
\[
 \int_{\mathbb T^{dr}}|H_a^{(r)}H_b^{(r)}|^{p/2}\,dm
 =\left(\int_{\mathbb T^d}|H_a^{(1)}H_b^{(1)}|^{p/2}\,dm\right)^r
 \le\kappa^r .
\]

Take the disjoint sets \(B_a^{(r)}=\{w:|H_a^{(r)}(w)|>|H_b^{(r)}(w)|\text{ for all }b\ne a\}\).
Outside \(B_a^{(r)}\), at least one \(b\ne a\) satisfies \(|H_b^{(r)}|\ge |H_a^{(r)}|\), so
\[
 \int_{(B_a^{(r)})^c}|H_a^{(r)}|^p\,dm
 \le\sum_{b\ne a}\int|H_a^{(r)}H_b^{(r)}|^{p/2}\,dm
 \le(|G|-1)\kappa^r\to0 \qquad\text{as }r\to\infty .
\]
For every \(r\) with \((|G|-1)\kappa^r<1\), we have
\(\|H_a^{(r)}1_{B_a^{(r)}}\|_{L^p(\mathbb T^{dr})}>0\). Define
\[
 v_a^{(r)}=\frac{H_a^{(r)}1_{B_a^{(r)}}}{\|H_a^{(r)}1_{B_a^{(r)}}\|_{L^p(\mathbb T^{dr})}}.
\]
By the reverse triangle inequality,
\[
 \|H_a^{(r)}-v_a^{(r)}\|_{L^p(\mathbb T^{dr})}
 \le2\|H_a^{(r)}1_{(B_a^{(r)})^c}\|_{L^p(\mathbb T^{dr})}
 \le2\big((|G|-1)\kappa^r\big)^{1/p}.
\]
Define \(S_0:L^p(G)\to L^p(\mathbb T^{dr})\) by
\[
 S_0 f=|G|^{-1/p}\sum_{a\in G} f(a)v_a^{(r)} .
\]
Integrating over the disjoint
sets \(B_a^{(r)}\) and using \(\|v_a^{(r)}\|_{L^p(\mathbb T^{dr})}=1\), we obtain
\[
 \|S_0f\|_{L^p(\mathbb T^{dr})}^p=|G|^{-1}\sum_a|f(a)|^p=\|f\|_{L^p(G)}^p .
\]
The triangle inequality and \(\|f\|_{L^1(G)}\le\|f\|_{L^p(G)}\) give
\[
 e_r:=\|S-S_0\|
 \le |G|^{1-1/p}\max_{a\in G}\|H_a^{(r)}-v_a^{(r)}\|_{L^p(\mathbb T^{dr})}
 \to0 \qquad\text{as }r\to\infty .
\]
We may therefore fix \(r\) sufficiently large that both \((|G|-1)\kappa^r<1\) and
\(e_r\le\varepsilon/(2+\varepsilon)\) hold.

\textbf{Step 2. Construct a left inverse for \(S\) with the required norm bound.}
By the definition of \(B_a^{(r)}\), \(|v_a^{(r)}(w)|>0\) for every \(w\in B_a^{(r)}\).
Define \(C:L^p(\mathbb T^{dr})\to L^p(G)\) by
\begin{equation}\label{eq:finite-recovery}
 (CF)(a)=|G|^{1/p}\int_{B_a^{(r)}} F\,\overline{v_a^{(r)}}|v_a^{(r)}|^{p-2}\,dm .
\end{equation}
For each \(a\in G\),
\[
 |G|^{-1/p}|(CF)(a)|
 \le \int_{B_a^{(r)}} |F|\,|v_a^{(r)}|^{p-1}\,dm
 \le \left(\int_{B_a^{(r)}}|F|^p\,dm\right)^{1/p}.
\]
The second inequality is an equality for \(p=1\) and follows from
H\"older's inequality and \(\|v_a^{(r)}\|_{L^p(\mathbb T^{dr})}=1\) for \(p>1\).
Since the sets \(B_a^{(r)}\) are pairwise disjoint,
\[
 \|CF\|_{L^p(G)}^p=\frac{1}{|G|}\sum_{a\in G}|(CF)(a)|^p
 \le \sum_{a\in G}\int_{B_a^{(r)}}|F|^p\,dm
 \le \|F\|_{L^p(\mathbb T^{dr})}^p .
\]
Thus \(\|C\|\le1\).
On \(B_a^{(r)}\) we have \(S_0f=|G|^{-1/p}f(a)v_a^{(r)}\), so
\[
 (CS_0 f)(a)=f(a)\int_{B_a^{(r)}}|v_a^{(r)}|^p\,dm=f(a).
\]
Hence \(CS_0=I\) and \(\|CS-I\|=\|C(S-S_0)\|\le e_r<1\), so the Neumann series gives the inverse
of \(CS\). Set \(Q=(CS)^{-1}C:L^p(\mathbb T^{dr})\to L^p(G)\). Then
\begin{equation}\label{eq:finite-left-inverse}
 QS=I,\qquad
 \|S\|\,\|Q\|
 \le\frac{1+e_r}{1-e_r}\le1+\varepsilon.
\end{equation}

\textbf{Step 3. Identify the monomials corresponding to the characters in \(A\).}
For \(h\in G\), let \(L_h f(a)=f(a-h)\) and let
\(T_hF(w)=F((\omega_k^{-h_k}w_{k,j})_{k,j})\).
The definition \eqref{eq:character-kernel} of \(H_a^{(r)}\) gives \(T_hH_a^{(r)}=H_{a+h}^{(r)}\).
The inequalities defining \(B_a^{(r)}\) then give
\[
 \begin{aligned}
 (T_h1_{B_a^{(r)}})(w)=1
 &\iff |(T_hH_a^{(r)})(w)|>|(T_hH_b^{(r)})(w)|
       \quad\text{for all }b\ne a\\
 &\iff |H_{a+h}^{(r)}(w)|>|H_c^{(r)}(w)|
       \quad\text{for all }c\ne a+h\\
 &\iff w\in B_{a+h}^{(r)}.
 \end{aligned}
\]
Thus \(T_h1_{B_a^{(r)}}=1_{B_{a+h}^{(r)}}\) and
\[
 T_hv_a^{(r)}
 =\frac{T_h(H_a^{(r)}1_{B_a^{(r)}})}
        {\|T_h(H_a^{(r)}1_{B_a^{(r)}})\|_{L^p(\mathbb T^{dr})}}
 =\frac{H_{a+h}^{(r)}1_{B_{a+h}^{(r)}}}
        {\|H_{a+h}^{(r)}1_{B_{a+h}^{(r)}}\|_{L^p(\mathbb T^{dr})}}
 =v_{a+h}^{(r)}.
\]
Substitution in \eqref{eq:finite-embedding} and a change of variables in \eqref{eq:finite-recovery} yield
\[
 T_h S=S L_h,\qquad C T_h=L_h C.
\]
Consequently \(CS\) and its inverse commute with every \(L_h\), and
\(Q T_h=L_h Q\).

For \(\nu=(\nu_1,\ldots,\nu_d)\in G\), the corresponding character is
\(\chi_\nu(a)=\prod_{k=1}^d\omega_k^{n_k m_k}\), where
\(\nu_k=[n_k]_{N_k}\) and \(a_k=[m_k]_{N_k}\).
Define \(\pi:\{0,1\}^{dr}\to G\) by
\[
 \pi(\alpha)_k
 =\left[\sum_{j=1}^r\alpha_{k,j}\right]_{N_k},
 \qquad 1\le k\le d .
\]
Set
\[
 D=\{\alpha\in\{0,1\}^{dr}:\chi_{\pi(\alpha)}\in A\},\qquad
 E_D=\operatorname{span}\{w^\alpha:\alpha\in D\}.
\]
For each \(h\in G\), \(L_h\chi_\nu=\overline{\chi_\nu(h)}\chi_\nu\) and
\(T_hw^\alpha=\overline{\chi_{\pi(\alpha)}(h)}w^\alpha\).
Since each \(H_a^{(r)}\) has degree at most one in each variable, we may write
\[
 (S\chi_\nu)(w)=\sum_{\alpha\in\{0,1\}^{dr}}c_\alpha w^\alpha .
\]
Applying \(T_h S=S L_h\) to \(\chi_\nu\) gives
\[
 \sum_{\alpha}c_\alpha\,\overline{\chi_{\pi(\alpha)}(h)}\,w^\alpha
 =T_h(S\chi_\nu)(w)
 =\overline{\chi_\nu(h)}\,(S\chi_\nu)(w)
 =\sum_{\alpha}c_\alpha\,\overline{\chi_\nu(h)}\,w^\alpha .
\]
Comparing the coefficients of \(w^\alpha\) and cancelling \(c_\alpha\) when
it is nonzero, we obtain
\[
 c_\alpha\ne0\quad\Longrightarrow\quad
 \chi_{\pi(\alpha)}(h)=\chi_\nu(h)\quad\text{for every }h\in G .
\]
Thus \(\pi(\alpha)=\nu\) whenever \(c_\alpha\ne0\). In particular, if
\(\chi_\nu\in A\), every such \(\alpha\) belongs to \(D\), so \(S\chi_\nu\in E_D\).
Consequently \(S(L^p_A(G))=\operatorname{span}\{S\chi:\chi\in A\}\subset E_D\).

Now fix \(\alpha\in D\) and put \(f=Q(w^\alpha)\). For all \(a,h\in G\),
\[
 f(a-h)=(L_h f)(a)=(Q(T_h w^\alpha))(a)
 =\overline{\chi_{\pi(\alpha)}(h)}\,f(a).
\]
Taking \(h=a\) gives \(f(0)=\overline{\chi_{\pi(\alpha)}(a)}\,f(a)\).
Since \(|\chi_{\pi(\alpha)}(a)|=1\), we obtain
\[
 Q(w^\alpha)(a)=f(0)\,\chi_{\pi(\alpha)}(a)\qquad(a\in G).
\]
For \(\alpha\in D\) we have \(\chi_{\pi(\alpha)}\in A\), and therefore
\(Q(w^\alpha)\in L^p_A(G)\). The monomials \(w^\alpha\) with \(\alpha\in D\)
span \(E_D\), hence \(Q(E_D)\subset L^p_A(G)\) by linearity.

\textbf{Step 4. Embed \(E_D\) isometrically into \(H^p(\mathbb T^{b+m})\).}
Put \(b=dr\). To write \(\alpha\) as a column vector in \(\mathbb Z^b\),
index its entries and the corresponding variables by
\[
 i=(k-1)r+j,\qquad \alpha_i=\alpha_{k,j},\qquad w_i=w_{k,j}
 \quad(1\le k\le d,\ 1\le j\le r).
\]
With \(m=b+|\{0,1\}^b\setminus D|\), Lemma~\ref{lem:cube-affine} gives
\(\Phi(\alpha)=c+B\alpha\in\mathbb Z^m\) with
\[
 \alpha\in D\quad\Longleftrightarrow\quad
 \alpha\ge0\ \text{and}\ \Phi(\alpha)\ge0\qquad(\alpha\in\mathbb Z^b).
\]

Write \(w=(w_1,\ldots,w_b)\in\mathbb T^b\), \(t=(t_1,\ldots,t_m)\in\mathbb T^m\),
and let \(B_i\in\mathbb Z^m\) be the \(i\)th column of \(B\).
For \(\beta\in\mathbb Z^m\), set \(t^\beta=\prod_{\ell=1}^m t_\ell^{\beta_\ell}\).
Define \(U:L^p(\mathbb T^b)\to L^p(\mathbb T^b\times\mathbb T^m)\) by
\[
 (Uf)(w,t)=t^c f(w_1t^{B_1},\ldots,w_bt^{B_b}),
\]
and \(V:L^p(\mathbb T^b\times\mathbb T^m)\to L^p(\mathbb T^b)\) by
\[
 (VF)(w)=\int_{\mathbb T^m}t^{-c}
 F(w_1t^{-B_1},\ldots,w_bt^{-B_b},t)\,dm(t).
\]
Integrating first in \(w\), we obtain
\[
 \|Uf\|_{L^p(\mathbb T^b\times\mathbb T^m)}^p
 =\int_{\mathbb T^m}\int_{\mathbb T^b}
 |f(w_1t^{B_1},\ldots,w_bt^{B_b})|^p\,dm(w)\,dm(t)
 =\|f\|_{L^p(\mathbb T^b)}^p .
\]
For \(V\), Jensen's inequality gives
\[
 |(VF)(w)|^p\le\int_{\mathbb T^m}
 |F(w_1t^{-B_1},\ldots,w_bt^{-B_b},t)|^p\,dm(t)
\]
for almost every \(w\). Integrating in \(w\) gives
\[
 \|VF\|_{L^p(\mathbb T^b)}^p
 \le\int_{\mathbb T^m}\int_{\mathbb T^b}|F(u,t)|^p\,dm(u)\,dm(t)
 =\|F\|_{L^p(\mathbb T^b\times\mathbb T^m)}^p .
\]
Moreover, for \(f\in L^p(\mathbb T^b)\),
\[
 (VUf)(w)=\int_{\mathbb T^m}t^{-c}t^c f(w)\,dm(t)=f(w)
 \qquad\text{for almost every }w\in\mathbb T^b.
\]
Thus \(VU=I\). On monomials,
\begin{equation}\label{eq:affine-monomials}
 U(w^\alpha)=w^\alpha t^{\Phi(\alpha)},\qquad
 V(w^\alpha t^\beta)=
 \begin{cases}w^\alpha,&\beta=\Phi(\alpha),\\0,&\beta\ne\Phi(\alpha).\end{cases}
\end{equation}
Hence, for \(F\in L^p(\mathbb T^b\times\mathbb T^m)\) and
\((\alpha,\beta)\in\mathbb Z^b\times\mathbb Z^m\),
\begin{equation}\label{eq:affine-fourier-projection}
 \widehat{UVF}(\alpha,\beta)=
 \begin{cases}
 \widehat F(\alpha,\beta),&\beta=\Phi(\alpha),\\
 0,&\beta\ne\Phi(\alpha).
 \end{cases}
\end{equation}
This follows for trigonometric polynomials by linearity and extends to
\(L^p\) by density, since \(UV\) and each Fourier coefficient functional
are bounded.
If \(\alpha\in D\), then \(\alpha\ge0\) and \(\Phi(\alpha)\ge0\), so
\(U(w^\alpha)=w^\alpha t^{\Phi(\alpha)}\) is an analytic monomial.

Now let \(w^\alpha t^\beta\) be an analytic monomial, so that
\(\alpha\ge0\) and \(\beta\ge0\). If \(\beta\ne\Phi(\alpha)\), then
\(V(w^\alpha t^\beta)=0\). If \(\beta=\Phi(\alpha)\), then
\(\Phi(\alpha)=\beta\ge0\). Together with \(\alpha\ge0\), this gives
\(\alpha\in D\) by \eqref{eq:cube-affine-criterion}, and hence
\(V(w^\alpha t^\beta)=w^\alpha\in E_D\).

\textbf{Step 5. Pass to the infinite torus and construct \(J\) and \(R\).}
Define
\[
 \begin{aligned}
 &\iota:L^p(\mathbb T^b\times\mathbb T^m)\longrightarrow L^p(\mathbb T^\infty),
 &&(\iota g)(w,t,z)=g(w,t),\\
 &P:L^p(\mathbb T^\infty)\longrightarrow L^p(\mathbb T^b\times\mathbb T^m),
 &&(PF)(w,t)=\int_{\mathbb T^\infty} F(w,t,z)\,dz.
 \end{aligned}
\]
Then
\(\|\iota g\|_{L^p(\mathbb T^\infty)}=\|g\|_{L^p(\mathbb T^b\times\mathbb T^m)}\)
and \(P\iota g=g\). Jensen's inequality in \(z\) gives
\(\|PF\|_{L^p(\mathbb T^b\times\mathbb T^m)}\le\|F\|_{L^p(\mathbb T^\infty)}\).

By \eqref{eq:affine-monomials}, for \(\alpha\in D\), \((\iota U)(w^\alpha)\) is the analytic monomial
\(w^\alpha t^{\Phi(\alpha)}\) on \(\mathbb T^\infty\). Hence
\((\iota U)(E_D)\subset H^p(\mathbb T^\infty)\).
For an analytic monomial \(w^\alpha t^\beta z^\gamma\),
\[
 P(w^\alpha t^\beta z^\gamma)=
 \begin{cases}
 w^\alpha t^\beta,&\gamma=0,\\
 0,&\gamma\ne0,
 \end{cases}
\]
because every nonconstant monomial in \(z\) has integral zero on \(\mathbb T^\infty\).
When \(\gamma=0\), \eqref{eq:cube-affine-criterion} and \eqref{eq:affine-monomials} give \(V(w^\alpha t^\beta)\in E_D\).
Since \(E_D\) is finite dimensional, it is closed.
Density of analytic polynomials and boundedness of \(VP\) now give
\((VP)(H^p(\mathbb T^\infty))\subset E_D\).

Set \(J=\iota US\) and \(R=QVP\). By \eqref{eq:finite-left-inverse},
\[
 RJ=QV(P\iota)US=Q(VU)S=I
\]
and
\[
 \|J\|\,\|R\|\le\|S\|\,\|Q\|\le 1+\varepsilon .
\]
The inclusions for \(S,\iota U\) give \(J(L^p_A(G))\subset H^p(\mathbb T^\infty)\),
and those for \(VP,Q\) give \(R(H^p(\mathbb T^\infty))\subset L^p_A(G)\).
In particular \(JR\) restricted to \(H^p(\mathbb T^\infty)\) is a projection
onto \(J(L^p_A(G))\).
\end{proof}

For \(A\ne\varnothing\), apply \eqref{eq:transfer} to
the restrictions of \(J\) and \(R\) in Theorem~\ref{thm:retract}. Using \eqref{eq:retraction-bounds} and letting
\(\varepsilon\downarrow0\) gives
\begin{equation}\label{eq:spectrum-lust}
 \lust(L_A^p(G))\le\lust(H^p(\mathbb T^\infty)).
\end{equation}
The case \(A=\varnothing\) follows from \(\lust(\{0\})=1\).

\subsection{A uniform \texorpdfstring{\(L^1(\mathbb T^b)\)}{L\string^1(T\string^b)} bound for partial sums in binary order}

For each positive integer \(b\), order the multi-indices in \(\{0,1\}^b\) by
\begin{equation}\label{eq:binary-order}
 \alpha\preceq\beta
 \quad\Longleftrightarrow\quad
 \sum_{j=1}^b 2^{j-1}\alpha_j\le \sum_{j=1}^b 2^{j-1}\beta_j .
\end{equation}

\begin{proposition}\label{prop:prefix}
Let \(D\) be a nonempty subset of \(\{0,1\}^b\). In the order \(\preceq\),
the monomials \(z^\alpha\), \(\alpha\in D\), form a basis of their span in
\(L^1(\mathbb T^b)\) with basis constant bounded independently of \(b\) and \(D\).
More precisely, there is an absolute constant \(C\) such that for every
\(f(z)=\sum_{\alpha\in D}a_\alpha z^\alpha\) and every \(\beta\in D\),
\begin{equation}\label{eq:binary-prefix-bound}
 \left\|\sum_{\substack{\alpha\in D\\\alpha\preceq\beta}}
 a_\alpha z^\alpha\right\|_{L^1(\mathbb T^b)}
 \le C\|f\|_{L^1(\mathbb T^b)} .
\end{equation}
\end{proposition}
\begin{proof}
Fix \(\beta\in D\), and denote the partial sum through \(\beta\) by
\(P_\beta f=\sum_{\alpha\in D,\,\alpha\preceq\beta}a_\alpha z^\alpha\).
Define an isometry on \(L^1(\mathbb T^b)\) by
\[
 (T_\beta f)(z)=z^\beta f(z_1^{1-2\beta_1},\ldots,z_b^{1-2\beta_b}).
\]
For \(\alpha\in D\), the exponent of \(z_j\) in \(T_\beta(z^\alpha)\) is
\(\beta_j+(1-2\beta_j)\alpha_j=|\alpha_j-\beta_j|\), since
\(\alpha_j,\beta_j\in\{0,1\}\). Thus, with \(\gamma(\alpha)=(|\alpha_j-\beta_j|)_{j=1}^b\),
\[
 g:=T_\beta f=\sum_{\alpha\in D}a_\alpha z^{\gamma(\alpha)}.
\]
Since \(\gamma(\beta)=0\), we have \(T_\beta(a_\beta z^\beta)=a_\beta\).
All other indices satisfy \(\gamma(\alpha)\ne0\), so \(a_\beta\) is the constant term of \(g\).
For \(\alpha\ne\beta\), let \(h(\alpha)=\max\{j:\alpha_j\ne\beta_j\}\).
Since \(2^{h(\alpha)-1}>\sum_{j<h(\alpha)}2^{j-1}\), \eqref{eq:binary-order} gives
\begin{equation}\label{eq:binary-order-criterion}
 \alpha\preceq\beta
 \quad\Longleftrightarrow\quad \alpha_{h(\alpha)}<\beta_{h(\alpha)}
 \quad\Longleftrightarrow\quad \beta_{h(\alpha)}=1 .
\end{equation}

For \(u\in L^1(\mathbb T^b)\), its conditional expectation \(C_j u\)
given the first \(j\) coordinates is
\[
 (C_j u)(z)=\int_{\mathbb T^{b-j}}
 u(z_1,\ldots,z_j,w_{j+1},\ldots,w_b)\,dm(w)
 \qquad (0\le j<b),
\]
with \(C_bu=u\). On a monomial \(z^\nu\), with \(\nu\in\{0,1\}^b\), this gives
\[
 C_j(z^\nu)
 =\prod_{i=1}^j z_i^{\nu_i}
  \prod_{i=j+1}^b\int_{\mathbb T}w_i^{\nu_i}\,dm(w_i)
 =\begin{cases}
 z^\nu,&\nu_i=0\text{ for all }i>j,\\
 0,&\text{otherwise}.
 \end{cases}
\]
Applying this formula to \(g=\sum_{\alpha\in D}a_\alpha z^{\gamma(\alpha)}\), we obtain
\[
\begin{aligned}
 &C_jg=a_\beta+\sum_{\substack{\alpha\in D\setminus\{\beta\}\\h(\alpha)\le j}}
       a_\alpha z^{\gamma(\alpha)}
       &&(0\le j\le b),\\
 &d_j:=C_jg-C_{j-1}g
 =\sum_{\substack{\alpha\in D\setminus\{\beta\}\\h(\alpha)=j}}
       a_\alpha z^{\gamma(\alpha)}
       &&(1\le j\le b).
\end{aligned}
\]
Since \(C_0g=a_\beta\), \eqref{eq:binary-order-criterion} yields
\begin{equation}\label{eq:binary-partial-sum}
 T_\beta(P_\beta f)
 =a_\beta+\sum_{\substack{\alpha\in D\setminus\{\beta\}\\\beta_{h(\alpha)}=1}}
       a_\alpha z^{\gamma(\alpha)}
 =C_0g+\sum_{\substack{1\le j\le b\\\beta_j=1}}d_j .
\end{equation}

For \(0\le k\le b\), set
\[
 M_k=C_kg-C_0g,\qquad N_k=\sum_{j=1}^k\beta_j d_j .
\]
Both are martingales starting at zero, with differences \(d_j\) and
\(\beta_jd_j\), respectively. Since \(|N_b|\le\max_{0\le k\le b}|N_k|\),
the square-function estimate \eqref{eq:martingale-square},
applied first to \(N\) and then to \(M\), gives
\[
 \|N_b\|_{L^1}
 \lesssim \left\|\left(\sum_{j=1}^b|\beta_j d_j|^2\right)^{1/2}\right\|_{L^1}
 \le \left\|\left(\sum_{j=1}^b|d_j|^2\right)^{1/2}\right\|_{L^1}
 \lesssim \left\|\max_{0\le k\le b}|C_kg-C_0g|\right\|_{L^1}.
\]

Analyticity of \(g\) gives \(C_kg=g(z_1,\ldots,z_k,0,\ldots,0)\).
The submean inequality for \(|g|^{1/2}\), applied in each of the last
\(b-k\) variables, gives \(|C_kg|^{1/2}\le C_k(|g|^{1/2})\).
Applying the \(L^2\) maximal inequality \eqref{eq:martingale-max} to the
martingale \((C_k(|g|^{1/2}))_{k=0}^b\), we obtain
\[
 \left\|\max_{0\le k\le b}|C_kg|\right\|_{L^1}
 \le \left\|\max_{0\le k\le b}C_k(|g|^{1/2})\right\|_{L^2}^{\,2}
 \le 4\big\||g|^{1/2}\big\|_{L^2}^{\,2}=4\|g\|_{L^1}.
\]
Since \(\max_k|C_kg-C_0g|\le 2\max_k|C_kg|\), the preceding estimates give
\(\|N_b\|_{L^1}\lesssim\|g\|_{L^1}\). By \eqref{eq:binary-partial-sum}, \(T_\beta(P_\beta f)=C_0g+N_b\).
Together with \(|C_0g|\le\|g\|_{L^1}\), this gives
\[
 \|P_\beta f\|_{L^1}
 =\|C_0g+N_b\|_{L^1}
 \le |C_0g|+\|N_b\|_{L^1}
 \lesssim \|g\|_{L^1}=\|f\|_{L^1},
\]
with an absolute implied constant.
\end{proof}

\subsection{Uniform bounds for polynomial embeddings and their left inverses}\label{sec:symmetric}

For each \(n\ge2\) set \(d=2^n-1\) and \(K=2d\). Let \(E_n\subset H^1(\mathbb T^d)\)
be the span of the monomials
\[
 z^\alpha=z_1^{\alpha_1}\cdots z_d^{\alpha_d},\qquad
 z=(z_1,\ldots,z_d)\in\mathbb T^d,\quad
 \alpha\in\mathbb N_0^d,\quad |\alpha|\le K.
\]
Let \(r\ge2K\) be an even integer and \(a_0\) an integer with
\(0\le a_0\le r-K\). Their values are specified in the proof of Lemma~\ref{lem:row-bound}.
To define the linear embedding
\(I_n:E_n\to L^1(\mathbb T^{dr})\) on the monomials \(z^\alpha\) with
\(|\alpha|\le K\), write the coordinates of \(w\in\mathbb T^{dr}\) as
\[
 w=\begin{pmatrix}
 w_{1,1}&\cdots&w_{1,r}\\
 \vdots&&\vdots\\
 w_{d,1}&\cdots&w_{d,r}
 \end{pmatrix},\qquad
 w_i=(w_{i,1},\ldots,w_{i,r})\in\mathbb T^r
 \quad(1\le i\le d).
\]
For a tuple of variables, \(e_k\) denotes the elementary symmetric polynomial
of degree \(k\), with \(e_0=1\). In particular, for \(1\le i\le d\),
\[
 e_{a_0+\alpha_i}(w_i)
 =\sum_{\substack{J\subseteq\{1,\ldots,r\}\\ |J|=a_0+\alpha_i}}
     \prod_{j\in J}w_{i,j}.
\]
Define \(I_n\) on monomials by
\begin{subequations}
\begin{equation}\label{eq:matrix-embedding}
 (I_n z^\alpha)(w)=A_{n,r}^{-d}\prod_{i=1}^d e_{a_0+\alpha_i}(w_i),
 \qquad |\alpha|\le K,
\end{equation}
where \(\xi=e^{2\pi i/r}\) and
\[
 A_{n,r}=r^{-r}\sum_{x\in\{0,\ldots,r-1\}^r}
             |e_{a_0}(\xi^{x_1},\ldots,\xi^{x_r})|.
\]
For \(\beta\in\{0,1\}^{d\times r}\) put \(w^\beta=\prod_{i,j}w_{i,j}^{\beta_{ij}}\).
Define \(\nu:\{0,1\}^{d\times r}\to\mathbb Z^d\) by
\[
 \nu_i(\beta)=\sum_{j=1}^r\beta_{ij}-a_0,\qquad 1\le i\le d .
\]
If \(w^\beta\) occurs in \(I_n z^\alpha\), then
\(\sum_{j=1}^r\beta_{ij}=a_0+\alpha_i\) for every \(i\), so \(\nu(\beta)=\alpha\).
Set
\begin{equation}\label{eq:matrix-indices}
 \mathcal D_n=\{\beta\in\{0,1\}^{d\times r}:\nu_i(\beta)\ge0\ (1\le i\le d),\
                   |\nu(\beta)|\le K\}.
\end{equation}
Let \(W_n=\operatorname{span}\{w^\beta:\beta\in\mathcal D_n\}
\subset L^1(\mathbb T^{dr})\). Then \(I_n(E_n)\subseteq W_n\).
Define \(M_n:W_n\to E_n\) by
\begin{equation}\label{eq:matrix-left-inverse}
 M_nw^\beta=t_\beta z^{\nu(\beta)},\qquad
 t_\beta=\frac{A_{n,r}^d}{\prod_{i=1}^d\binom r{a_0+\nu_i(\beta)}}.
\end{equation}
\end{subequations}
Let \(\mathrm{Sym}(r)\) be the permutation group on \(\{1,\ldots,r\}\).
By the definition of \(\nu\),
\begin{equation}\label{eq:matrix-count}
 \#\{\beta\in\{0,1\}^{d\times r}:\nu(\beta)=\alpha\}
 =\prod_{i=1}^d\binom{r}{a_0+\alpha_i}
 \qquad(\alpha\in\mathbb N_0^d,\ |\alpha|\le K).
\end{equation}
Hence \(M_n I_n=I_{E_n}\), and for \(F\in W_n\),
\begin{equation}\label{eq:row-permutations}
 (I_nM_n F)(w)=(r!)^{-d}\!
 \sum_{\pi\in(\mathrm{Sym}(r))^d}F\big((w_{i,\pi_i(j)})_{i,j}\big).
\end{equation}

\begin{lemma}\label{lem:row-bound}
The parameters \(r,a_0\) can be chosen explicitly as functions of \(n\)
so that the operators in \eqref{eq:matrix-embedding} and \eqref{eq:matrix-left-inverse} satisfy
\(\sup_{n\ge2}(\|I_n:E_n\to W_n\|+\|M_n:W_n\to E_n\|)<\infty\).
\end{lemma}
\begin{proof}
Let \(\mathcal P_K=\{h(v)=\sum_{k=0}^K b_kv^k:b_0,\ldots,b_K\in\mathbb C\}\),
with \(v\in\mathbb T\).
For \(u=(u_1,\ldots,u_r)\in\mathbb T^r\), define
\begin{equation}\label{eq:row-operator}
 U_r:\mathcal P_K\longrightarrow L^1(\mathbb T^r),\qquad
 (U_rh)(u)=(\pi/4)^r\sum_{k=0}^K b_k e_{a_0+k}(u).
\end{equation}
Write \(u/v=(u_1/v,\ldots,u_r/v)\) and set
\(Q_r(u)=(\pi/4)^r\prod_{j=1}^r(1+u_j)\). Expanding this product gives
\(Q_r(u/v)=(\pi/4)^r\sum_{\ell=0}^r e_\ell(u)v^{-\ell}\). Integrating in \(v\)
therefore gives
\begin{equation}\label{eq:row-integral}
 (U_rh)(u)=\int_{\mathbb T}h(v)v^{a_0}Q_r(u/v)\,dm(v).
\end{equation}

For \(u\in\mathbb T^r\), let \(Z(u)=(\sum_j u_j)/|\sum_j u_j|\) when
\(\sum_j u_j\ne0\), and let \(Z(u)=u_1\) otherwise. Taking \(h\equiv1\)
in \eqref{eq:row-operator} gives \(U_r1=(\pi/4)^r e_{a_0}\).
Thus \((U_r1)(vu)=v^{a_0}(U_r1)(u)\) by homogeneity of \(e_{a_0}\),
where \(vu=(vu_1,\ldots,vu_r)\). Since \(Z(vu)=vZ(u)\), averaging over \(v\) gives
\begin{equation}\label{eq:row-product-norm}
 \begin{aligned}
 \|(U_r1)\,(h\circ Z)\|_{L^1(\mathbb T^r)}
 &=\int_{\mathbb T^r}|(U_r1)(u)|\int_{\mathbb T}|h(vZ(u))|\,dm(v)\,dm(u)\\
 &=\|U_r1\|_{L^1(\mathbb T^r)}\|h\|_{L^1(\mathbb T)}.
 \end{aligned}
\end{equation}
Applying \eqref{eq:row-integral} to \(h\) and to \(1\) gives
\[
 (U_rh)(u)-(U_r1)(u)h(Z(u))
 =\int_{\mathbb T}\bigl(h(v)-h(Z(u))\bigr)v^{a_0}Q_r(u/v)\,dm(v).
\]
Hence
\[
 \begin{aligned}
 \|U_rh-(U_r1)\,(h\circ Z)\|_{L^1(\mathbb T^r)}
 &\le\int_{\mathbb T}\int_{\mathbb T^r}|h(v)-h(Z(u))|\,|Q_r(u/v)|\,dm(u)\,dm(v)\\
 &=\int_{\mathbb T^r}|Q_r(x)|
      \int_{\mathbb T}|h(v)-h(vZ(x))|\,dm(v)\,dm(x).
 \end{aligned}
\]
For each \(x\in\mathbb T^r\), write \(Z(x)=e^{it}\) with \(t=\arg Z(x)\in(-\pi,\pi]\).
The fundamental theorem of calculus and Bernstein's inequality \cite{Arestov} give
\[
 \int_{\mathbb T}|h(v)-h(vZ(x))|\,dm(v)
 \le |t|\left\|\frac{d}{d\theta}h(e^{i\theta})\right\|_{L^1(\mathbb T)}
 \le K|t|\,\|h\|_{L^1(\mathbb T)}.
\]
Inserting this bound for the inner integral gives
\begin{equation}\label{eq:row-error}
 \|U_rh-(U_r1)\,(h\circ Z)\|_{L^1(\mathbb T^r)}
 \le K\|h\|_{L^1(\mathbb T)}
       \int_{\mathbb T^r}|\arg Z(x)|\,|Q_r(x)|\,dm(x).
\end{equation}

To estimate the integral on the right, write \(x_\ell=e^{iX_\ell}\) with
\(X_\ell=\arg x_\ell\in(-\pi,\pi]\). For \(1\le j\le r\), set
\(d\mu_j(x)=\prod_{\ell=1}^j\frac{\cos(X_\ell/2)}4\,dX_\ell\)
on \(\mathbb T^j\), so that \(d\mu_r(x)=|Q_r(x)|\,dm(x)\).
Define \(T_j:\mathbb T^j\to(-\pi,\pi]\) by
\[
 T_j(x_1,\ldots,x_j)=
 \begin{cases}
 \arg(x_1+\cdots+x_j), & x_1+\cdots+x_j\ne0,\\
 \arg x_1, & x_1+\cdots+x_j=0.
 \end{cases}
\]
Since \(T_r(x)=\arg Z(x)\), the integral to be estimated is
\(\int_{\mathbb T^r}|T_r(x)|\,d\mu_r(x)\).
Integration in one coordinate gives
\[
 \int_{\mathbb T}\cos X_1\,d\mu_1=\frac13,\qquad
 \int_{\mathbb T}\sin X_1\,d\mu_1=0,\qquad
 \int_{\mathbb T}\sin^2 X_1\,d\mu_1=\frac{8}{15}.
\]
Write \(j^{-1}\sum_{\ell=1}^j x_\ell=C_j+iS_j\), where
\[
 C_j=\frac1j\sum_{\ell=1}^j\cos X_\ell,\qquad
 S_j=\frac1j\sum_{\ell=1}^j\sin X_\ell.
\]
For \(k\ne\ell\), the integral of \(\sin X_k\sin X_\ell\) vanishes by the product
formula for \(d\mu_j\). Thus
\[
 \int_{\mathbb T^j}S_j^2\,d\mu_j
 =\frac1{j^2}\sum_{\ell=1}^j\int_{\mathbb T^j}\sin^2 X_\ell\,d\mu_j
 =\frac{8}{15j}.
\]
Since \(\int_{\mathbb T^j}C_j\,d\mu_j=1/3\) and \(-1\le-\cos X_\ell\le1\),
Hoeffding's inequality \cite[Theorem~2]{Hoeffding} applied to the product measure
\(\mu_j\) gives
\[
 \int_{\{C_j<3/10\}}d\mu_j
 =\int_{\{1/3-C_j>1/30\}}d\mu_j
 \le \exp\!\left(-\frac{2j(1/30)^2}{2^2}\right)=e^{-j/1800}.
\]
On \(\{C_j\ge3/10\}\), the real part is positive, so
\[
 |T_j(x)|=\left|\arctan\frac{S_j}{C_j}\right|
 \le\frac{|S_j|}{C_j}\le\frac{10}{3}|S_j|.
\]
On \(\{C_j<3/10\}\), we use \(|T_j(x)|\le\pi\). Integrating over the two sets gives
\[
 \begin{aligned}
 \int_{\mathbb T^j}T_j(x)^2\,d\mu_j(x)
 &=\int_{\{C_j\ge3/10\}}T_j(x)^2\,d\mu_j(x)
   +\int_{\{C_j<3/10\}}T_j(x)^2\,d\mu_j(x)\\
 &\le \frac{100}{9}\int_{\mathbb T^j}S_j^2\,d\mu_j
       +\pi^2\int_{\{C_j<3/10\}}d\mu_j
 \le \frac{160}{27j}+\pi^2e^{-j/1800}.
 \end{aligned}
\]
For \(r\ge10^6\), taking \(j=r\) and \(j=r-1\) gives
\begin{equation}\label{eq:angle-second-moments}
 \int_{\mathbb T^r}T_r(x)^2\,d\mu_r(x)\le\frac9r,\qquad
 r\int_{\mathbb T^{r-1}}T_{r-1}(y)^2\,d\mu_{r-1}(y)<6.
\end{equation}
The first bound and the Cauchy--Schwarz inequality give
\begin{equation}\label{eq:angle-first-moment}
 \int_{\mathbb T^r}|T_r(x)|\,d\mu_r(x)\le\frac{3}{\sqrt r}.
\end{equation}

To obtain the lower bound for \(\|U_r1\|_{L^1(\mathbb T^r)}\), write
\(\operatorname{sgn}z=z/|z|\) for \(z\ne0\) and \(\operatorname{sgn}0=0\).
Since \(|\operatorname{sgn}z|\le1\), \eqref{eq:row-integral} with \(h=1\) and the substitution
\(u=vx\) give
\[
\begin{aligned}
 \|U_r1\|_{L^1(\mathbb T^r)}
 &\ge\operatorname{Re}\int_{\mathbb T^r}(U_r1)(u)\,Z(u)^{-a_0}
       \overline{\operatorname{sgn}Q_r(u/Z(u))}\,dm(u)\\
 &=\operatorname{Re}\int_{\mathbb T^r}
 Z(x)^{-a_0}\operatorname{sgn}Q_r(x)\,
 \overline{\operatorname{sgn}Q_r(x/Z(x))}\,d\mu_r(x).
\end{aligned}
\]
Let \(B=\{x\in\mathbb T^r:\max_{1\le j\le r}|X_j-T_r(x)|\ge\pi\}\).
For \(-\pi<t<\pi\), \(\operatorname{sgn}(1+e^{it})=e^{it/2}\).
On \(B^c\), we have \(-\pi<X_j-T_r(x)<\pi\) for every \(j\).
Since \(x_j=e^{iX_j}\) and \(x_j/Z(x)=e^{i(X_j-T_r(x))}\), the product formula
for \(Q_r\) gives, for \(\mu_r\)-almost every \(x\in B^c\),
\[
 \operatorname{sgn}Q_r(x)=\prod_{j=1}^r e^{iX_j/2},\qquad
 \operatorname{sgn}Q_r(x/Z(x))=\prod_{j=1}^r e^{i(X_j-T_r(x))/2}.
\]
Multiplying these factors on \(B^c\) gives
\[
 \begin{aligned}
 Z(x)^{-a_0}\operatorname{sgn}Q_r(x)\,
 \overline{\operatorname{sgn}Q_r(x/Z(x))}
 &=e^{-ia_0T_r(x)}\prod_{j=1}^r e^{iX_j/2}e^{-i(X_j-T_r(x))/2}\\
 &=e^{-ia_0T_r(x)}e^{irT_r(x)/2}.
 \end{aligned}
\]
Choose \(a_0=r/2\). The product is then \(1\) for \(\mu_r\)-almost every \(x\in B^c\).
On \(B\) its modulus is at most \(1\), so its real part is at least \(-1\).
Integrating over \(B^c\) and \(B\) gives
\begin{equation}\label{eq:row-lower-exception}
 \|U_r1\|_{L^1(\mathbb T^r)}
 \ge \int_{B^c}d\mu_r-\int_B d\mu_r=1-2\int_B d\mu_r.
\end{equation}

We now bound \(\int_B d\mu_r\) from above. For \(x=(x_1,\ldots,x_r)\in\mathbb T^r\)
and \(1\le i\le r\), set
\(t_i(x)=T_{r-1}(x_1,\ldots,x_{i-1},x_{i+1},\ldots,x_r)\).
If \(|X_i-t_i(x)|<\pi\), the sector between these angles is convex, so
\[
 T_r(x)=\arg\!\left(e^{iX_i}+
        \left|\sum_{j\ne i}x_j\right|e^{it_i(x)}\right)
 \in \bigl[\min\{X_i,t_i(x)\},\,\max\{X_i,t_i(x)\}\bigr].
\]
Thus \(|X_i-T_r(x)|\le |X_i-t_i(x)|<\pi\). Taking complements gives
\[
 \{x:|X_i-T_r(x)|\ge\pi\}\subseteq\{x:|X_i-t_i(x)|\ge\pi\}.
\]
Since \(B=\bigcup_{i=1}^r\{x:|X_i-T_r(x)|\ge\pi\}\),
the preceding inclusion gives
\[
 \int_B d\mu_r
 \le\sum_{i=1}^r\int_{\{x:\,|X_i-t_i(x)|\ge\pi\}}d\mu_r(x).
\]
By Fubini's theorem,
\[
 \int_{\{x:\,|X_i-t_i(x)|\ge\pi\}}d\mu_r(x)
 =\int_{\mathbb T^{r-1}}\!
   \left(\int_{\substack{-\pi<X_i\le\pi\\ |X_i-t_i(x)|\ge\pi}}
   \frac{\cos(X_i/2)}4\,dX_i\right)\prod_{\ell\ne i}d\mu_1(x_\ell).
\]
The density is even, so the inner integral equals
\[
 \int_{\pi-|t_i(x)|}^{\pi}\frac{\cos(X_i/2)}4\,dX_i
 =\frac{1-\cos(t_i(x)/2)}2\le\frac{t_i(x)^2}{16}.
\]
Integrating over the remaining coordinates, summing over \(i\), and using \eqref{eq:angle-second-moments}, we obtain
\[
 \int_B d\mu_r
 \le\frac{r}{16}\int_{\mathbb T^{r-1}}T_{r-1}(x)^2\,d\mu_{r-1}(x)<\frac38.
\]
Substituting this estimate into \eqref{eq:row-lower-exception} gives
\[
 \|U_r1\|_{L^1(\mathbb T^r)}
 \ge 1-2\int_B d\mu_r>1-2\cdot\frac38=\frac14.
\]

Now set \(A_*=\|e_{a_0}\|_{L^1(\mathbb T^r)}=(4/\pi)^r\|U_r1\|_{L^1(\mathbb T^r)}\).
By \eqref{eq:row-product-norm} and the reverse triangle inequality, followed by
\eqref{eq:row-error} and \eqref{eq:angle-first-moment}, we have
\[
 \begin{aligned}
 \big|\|U_rh\|_{L^1(\mathbb T^r)}
 -\|U_r1\|_{L^1(\mathbb T^r)}\|h\|_{L^1(\mathbb T)}\big|
 &\le \|U_rh-(U_r1)(h\circ Z)\|_{L^1(\mathbb T^r)}\\
 &\le \frac{3K}{\sqrt r}\|h\|_{L^1(\mathbb T)}.
 \end{aligned}
\]
Dividing by \(\|U_r1\|_{L^1(\mathbb T^r)}>1/4\) and using \eqref{eq:row-operator}, we obtain
\begin{equation}\label{eq:one-variable-embedding}
 \left(1-\frac{12K}{\sqrt r}\right)\|h\|_{L^1(\mathbb T)}
 \le \left\|\sum_{k=0}^K b_k\,\frac{e_{a_0+k}}{A_*}\right\|_{L^1(\mathbb T^r)}
 \le \left(1+\frac{12K}{\sqrt r}\right)\|h\|_{L^1(\mathbb T)} .
\end{equation}

For \(f(z)=\sum_{|\alpha|\le K} b_\alpha z^\alpha\in E_n\) and
\(0\le j\le d\), define \(F_j\) on
\(\Omega_j=(\mathbb T^r)^j\times\mathbb T^{d-j}\) by
\[
 F_j(w_1,\ldots,w_j,z_{j+1},\ldots,z_d)
 =\sum_{|\alpha|\le K} b_\alpha
   \prod_{i=1}^j\frac{e_{a_0+\alpha_i}(w_i)}{A_*}
   \prod_{i=j+1}^d z_i^{\alpha_i}.
\]
Fix \(j\in\{1,\ldots,d\}\) and the coordinates \(w_1,\ldots,w_{j-1},z_{j+1},\ldots,z_d\).
Collecting the terms of \(F_{j-1}\) with \(\alpha_j=k\) gives the coefficient of \(z_j^k\),
\[
 c_k=\sum_{\substack{|\alpha|\le K\\ \alpha_j=k}} b_\alpha
       \prod_{i<j}\frac{e_{a_0+\alpha_i}(w_i)}{A_*}
       \prod_{i>j}z_i^{\alpha_i}
       \qquad (0\le k\le K).
\]
With these same coefficients, the definitions of \(F_{j-1}\) and \(F_j\) give
\[
 \begin{aligned}
 &F_{j-1}(w_1,\ldots,w_{j-1},z_j,\ldots,z_d)
   =\sum_{k=0}^K c_k z_j^k,\\
 &F_j(w_1,\ldots,w_j,z_{j+1},\ldots,z_d)
   =\sum_{k=0}^K c_k\,\frac{e_{a_0+k}(w_j)}{A_*}.
 \end{aligned}
\]
Applying \eqref{eq:one-variable-embedding} to \(h(v)=\sum_{k=0}^K c_kv^k\) therefore gives
\[
 \left(1-\frac{12K}{\sqrt r}\right)\int_{\mathbb T}|F_{j-1}|\,dm(z_j)
 \le \int_{\mathbb T^r}|F_j|\,dm(w_j)
 \le \left(1+\frac{12K}{\sqrt r}\right)\int_{\mathbb T}|F_{j-1}|\,dm(z_j).
\]
Integrating over \((w_1,\ldots,w_{j-1},z_{j+1},\ldots,z_d)\) in
\((\mathbb T^r)^{j-1}\times\mathbb T^{d-j}\) and applying Fubini's theorem yields
\[
 \left(1-\frac{12K}{\sqrt r}\right)\|F_{j-1}\|_{L^1(\Omega_{j-1})}
 \le \|F_j\|_{L^1(\Omega_j)}
 \le\left(1+\frac{12K}{\sqrt r}\right)\|F_{j-1}\|_{L^1(\Omega_{j-1})}.
\]
When \(12K/\sqrt r<1\), both factors are positive. Starting with \(F_0=f\)
and applying the two inequalities successively for \(j=1,\ldots,d\) gives
\begin{equation}\label{eq:iterated-row-embedding}
 \left(1-\frac{12K}{\sqrt r}\right)^d\|f\|_{L^1(\mathbb T^d)}
 \le \|F_d\|_{L^1(\mathbb T^{dr})}
 \le \left(1+\frac{12K}{\sqrt r}\right)^d\|f\|_{L^1(\mathbb T^d)}.
\end{equation}
To keep both factors bounded away from zero and infinity, choose
\begin{equation}\label{eq:block-parameters}
 r=(12d)^4,\qquad a_0=r/2.
\end{equation}
Then \(r\) is even and \(r\ge36^4>10^6\), since \(d\ge3\).
As \(K=2d\), we have \(12K/\sqrt r=1/(6d)\), so
\[
 \left(1-\frac{1}{6d}\right)^d\ge 1-\frac16=\frac56,\qquad
 \left(1+\frac{1}{6d}\right)^d\le\left(1-\frac{1}{6d}\right)^{-d}\le\frac65 .
\]

By \eqref{eq:matrix-embedding}, \(I_nf=(A_*/A_{n,r})^d F_d\).
To pass from bounds for \(\|F_d\|_{L^1(\mathbb T^{dr})}\) to bounds for \(\|I_nf\|_{L^1(\mathbb T^{dr})}\), we must therefore bound
\((A_{n,r}/A_*)^d\) above and away from zero independently of \(n\).
By definition,
\[
 \begin{aligned}
 &A_*=\int_{\mathbb T}\!\cdots\!\int_{\mathbb T}
       |e_{a_0}(w_1,\ldots,w_r)|\,dm(w_1)\cdots dm(w_r),\\
 &A_{n,r}=\frac{1}{r^r}\sum_{k_1=0}^{r-1}\cdots\sum_{k_r=0}^{r-1}
       |e_{a_0}(\xi^{k_1},\ldots,\xi^{k_r})|.
 \end{aligned}
\]
Fix \(w_1,\ldots,w_{r-1}\in\mathbb T\) and set \(w_r=z\).
In the defining sum for \(e_{a_0}\), separate the products that contain \(w_r\)
from those that omit it. This gives
\[
 \begin{aligned}
 p(z)&:=e_{a_0}(w_1,\ldots,w_{r-1},z)\\
 &=\underbrace{\sum_{\substack{J\subseteq\{1,\ldots,r-1\}\\ |J|=a_0}}
       \prod_{\ell\in J}w_\ell}_{a\,=\,e_{a_0}(w_1,\ldots,w_{r-1})}
 +z\,\underbrace{\sum_{\substack{J\subseteq\{1,\ldots,r-1\}\\ |J|=a_0-1}}
       \prod_{\ell\in J}w_\ell}_{b\,=\,e_{a_0-1}(w_1,\ldots,w_{r-1})}.
 \end{aligned}
\]
Set \(h=2\pi/r\), so that \(\xi=e^{ih}\).
Taylor's theorem with integral remainder gives \(|e^{is}-1-is|\le s^2/2\) for real \(s\).
For \(z\in\mathbb T\) and \(0\le u\le1\), applying this estimate at \(s=-uh\) and \(s=(1-u)h\)
and using \((1-u)(-uh)+u(1-u)h=0\) gives
\[
 \begin{aligned}
 \big|p(ze^{iuh})-(1-u)p(z)-u p(z\xi)\big|
 &=|b|\,\big|1-(1-u)e^{-iuh}-u e^{i(1-u)h}\big|\\
 &\le \tfrac12|b|h^2\big((1-u)u^2+u(1-u)^2\big)
 =\tfrac12|b|h^2u(1-u).
 \end{aligned}
\]
Write \(S(t)=r^{-1}\sum_{k=0}^{r-1}|p(e^{it}\xi^k)|\), so \(S(t+h)=S(t)\).
Averaging the last estimate over \(z=e^{it}\xi^k\) gives
\[
 S(t+uh)\le (1-u)S(t)+uS(t+h)+\tfrac12|b|h^2u(1-u)
 \le S(t)+\tfrac18|b|h^2.
\]
Taking \(t=0\) bounds \(S(uh)-S(0)\). For the reverse difference use
\(S(0)-S(uh)=S(uh+(1-u)h)-S(uh)\). Thus
\[
 |S(uh)-S(0)|\le |b|h^2/8,\qquad 0\le u\le1.
\]
Since
\[
 \int_0^1 S(uh)\,du
 =\frac1{rh}\sum_{k=0}^{r-1}\int_{kh}^{(k+1)h}|p(e^{it})|\,dt
 =\|p\|_{L^1(\mathbb T)}
\]
and \(|b|=\left|\int_{\mathbb T}p(z)\bar z\,dm(z)\right|\le\|p\|_{L^1(\mathbb T)}\),
integration in \(u\) yields
\[
 \left|\frac1r\sum_{k=0}^{r-1}|p(\xi^k)|-\|p\|_{L^1(\mathbb T)}\right|
 \le \frac{|b|h^2}{8}
 \le \frac{\pi^2}{2r^2}\|p\|_{L^1(\mathbb T)}.
\]
For \(0\le j\le r\), put
\[
 A_{n,j}=\frac{1}{r^j}\sum_{k_1,\ldots,k_j=0}^{r-1}
 \int_{\mathbb T^{r-j}}
 \big|e_{a_0}(\xi^{k_1},\ldots,\xi^{k_j},w_{j+1},\ldots,w_r)\big|
 \prod_{\ell=j+1}^r dm(w_\ell).
\]
Then \(A_{n,0}=A_*\). By symmetry of \(e_{a_0}\), the
one-variable estimate applies to
\[
 z\longmapsto e_{a_0}(\xi^{k_1},\ldots,\xi^{k_{j-1}},z,w_{j+1},\ldots,w_r).
\]
Summing in \(k_1,\ldots,k_{j-1}\) with weight \(r^{1-j}\) and integrating in
\(w_{j+1},\ldots,w_r\) gives
\[
 \left(1-\frac{\pi^2}{2r^2}\right)A_{n,j-1}
 \le A_{n,j}\le
 \left(1+\frac{\pi^2}{2r^2}\right)A_{n,j-1}
 \qquad (1\le j\le r).
\]
Iterating from \(j=1\) to \(r\) and taking \(d\)th powers yields
\[
 \left(1-\frac{\pi^2}{2r^2}\right)^{rd}
 \le \left(\frac{A_{n,r}}{A_*}\right)^d
 \le \left(1+\frac{\pi^2}{2r^2}\right)^{rd}.
\]
By \eqref{eq:block-parameters}, \(\pi^2d/(2r)<1/6\). Hence
\[
 \begin{aligned}
 &\left(1-\frac{\pi^2}{2r^2}\right)^{rd}
 \ge 1-\frac{\pi^2 d}{2r}>\frac56,\\
 &\left(1+\frac{\pi^2}{2r^2}\right)^{rd}
 \le\left(1-\frac{\pi^2}{2r^2}\right)^{-rd}<\frac65.
 \end{aligned}
\]
Substituting these bounds in \(I_nf=(A_*/A_{n,r})^d F_d\) and using
\eqref{eq:iterated-row-embedding} yields
\[
 \left(\frac56\right)^2\|f\|_{L^1(\mathbb T^d)}
 \le\|I_nf\|_{L^1(\mathbb T^{dr})}
 \le\left(\frac65\right)^2\|f\|_{L^1(\mathbb T^d)}.
\]
In particular,
\(\|I_n\|\le2\) and \(\|I_n^{-1}:I_n(E_n)\to E_n\|\le2\).
For \(F\in W_n\), \eqref{eq:row-permutations} gives
\[
 \|I_nM_n F\|_{L^1(\mathbb T^{dr})}
 \le (r!)^{-d}\!\sum_{\pi\in\mathrm{Sym}(r)^d}
 \big\|F((w_{i,\pi_i(j)})_{i,j})\big\|_{L^1(\mathbb T^{dr})}
 =\|F\|_{L^1(\mathbb T^{dr})}.
\]
It follows that \(\|M_n\|\le\|I_n^{-1}\|\,\|I_nM_n\|\le2\) for every \(n\ge2\).
\end{proof}

\section{Proof of Theorem~\ref*{thm:unc}}\label{sec:unc-proof}

\subsection{Fourier support outside the graph of a quadratic form}\label{sec:mixed}

We adapt the two-dimensional blocks and coefficient estimates of Komorowski and Tomczak-Jaegermann \cite[proof of Theorem~3.1]{KT}. To apply the result of Section~\ref{sec:finite}, we construct a finite character set \(A\) and a projection on \(L_A^p(G)\) whose range has the required blocks. The projection norm will depend only on \(p\).

Fix \(1<p<\infty\), \(p\ne2\), and a positive integer \(n\).
Let \(\mathbb F_3=\mathbb Z/3\mathbb Z\), \(V=\mathbb F_3^n\), and
\(\omega=e^{2\pi i/3}\). Then \(|V|=3^n\).
For \(u,v\in V\), write \(u\cdot v=\sum_{\ell=1}^n u_\ell v_\ell\in\mathbb F_3\).
The value of \(\omega^a\) is independent of the integer representing \(a\in\mathbb F_3\).
Set \(G=\mathbb F_3\times V^3\). Its characters are
\[
 \chi_{\tau,\xi}(s,x)=\omega^{\,s\tau+\sum_{j=1}^3 x_j\cdot\xi_j},
 \qquad \tau,s\in\mathbb F_3,\quad
 x=(x_1,x_2,x_3),\ \xi=(\xi_1,\xi_2,\xi_3)\in V^3 .
\]
Define \(q:V\to\mathbb F_3\) by \(q(v)=v\cdot v\), and
\(b:V^3\to\mathbb F_3\) by \(b(\xi)=-\sum_{j=1}^3q(\xi_j)\).
Set
\[
 A=\{\chi_{\tau,\xi}:\tau\in\mathbb F_3,\ \xi\in V^3,\ \tau\ne b(\xi)\}
 \subset\widehat G.
\]
With the norm inherited from \(L^p(G)\), the corresponding space is
\[
 E=L_A^p(G)=\left\{\sum_{\xi\in V^3}
       \sum_{\substack{\tau\in\mathbb F_3\\ \tau\ne b(\xi)}}
       c_{\tau,\xi}\chi_{\tau,\xi}:c_{\tau,\xi}\in\mathbb C\right\}.
\]
Write \(\widehat F(\tau,\xi)\) for the Fourier coefficient at
\(\chi_{\tau,\xi}\). The Fourier inversion formula gives
\[
 F\in E\quad\Longleftrightarrow\quad
 \widehat F(b(\xi),\xi)=0\quad\text{for every }\xi\in V^3 .
\]
On \(V\), the Fourier index \(\eta\in V\) denotes the character
\(x\mapsto\omega^{x\cdot\eta}\). On \(V^3\) we use the corresponding product of three characters.
Define \(T:L^p(V)\to L^p(V)\) and \(U_s:L^p(V^3)\to L^p(V^3)\) by
\[
 \widehat{Tf}(\eta)=\omega^{q(\eta)}\widehat f(\eta)\quad(\eta\in V),
 \qquad
 \widehat{U_sg}(\xi)=\omega^{\sum_{j=1}^3r_{sj}q(\xi_j)}\widehat g(\xi)
 \quad(\xi\in V^3),
\]
where
\[
 (r_{sj})_{0\le s\le2,\,1\le j\le3}=
 \begin{pmatrix}0&1&-1\\-1&0&1\\1&-1&0\end{pmatrix}.
\]
Write \(F_s(x)=F(s,x)\) and let \(\widehat F_s(\xi)\) denote its normalized
Fourier coefficient on \(V^3\). Since \(|G|=3|V|^3\), \eqref{eq:finite-fourier} gives
\[
 \widehat F(\tau,\xi)=\frac13\sum_{s=0}^2\omega^{-s\tau}
 \underbrace{\frac1{|V|^3}\sum_{x\in V^3}F_s(x)
       \omega^{-\sum_{j=1}^3 x_j\cdot\xi_j}}_{\widehat F_s(\xi)} .
\]
The identity \(r_{sj}=r_{0j}-s\) in \(\mathbb F_3\) gives
\(\sum_j r_{sj}q(\xi_j)=\sum_j r_{0j}q(\xi_j)+s b(\xi)\).
Taking the inverse multiplier therefore gives
\[
 \sum_s\widehat{U_s^{-1}F_s}(\xi)
 =\omega^{-\sum_j r_{0j}q(\xi_j)}
 \sum_s\omega^{-s b(\xi)}\widehat F_s(\xi)
 =3\omega^{-\sum_j r_{0j}q(\xi_j)}\widehat F(b(\xi),\xi).
\]
Consequently,
\begin{equation}\label{eq:quadratic-fibers}
 F\in E\quad\Longleftrightarrow\quad \sum_{s=0}^2 U_s^{-1}F_s=0.
\end{equation}

\subsection{A projection and its three norms}
For \(1\le j\le n\), define \(e_j\in V=\mathbb F_3^n\) by
\[
 (e_j)_k=\begin{cases}1,&k=j,\\0,&k\ne j,\end{cases}
 \qquad 1\le k\le n.
\]
Set \(I=\{e_1,\ldots,e_n\}\) and let \(P:L^p(V)\to L^p(V)\) be given by \(Pf=1_I f\).
For \(v,\xi\in V\), the Fourier coefficient of \(1_{\{v\}}\) in \eqref{eq:finite-fourier} is
\[
 \widehat{1_{\{v\}}}(\xi)
 =\frac{1}{|V|}\sum_{y\in V}1_{\{v\}}(y)\omega^{-y\cdot\xi}
 =\frac{\omega^{-v\cdot\xi}}{|V|}.
\]
By the definition of \(T\), we have
\(\widehat{T^r f}(\xi)=\omega^{r q(\xi)}\widehat f(\xi)\) for \(r=\pm1\).
Applying the Fourier inversion formula to \(T^r1_{\{v\}}\), we obtain
\[
 (T^r1_{\{v\}})(x)
 =\sum_{\xi\in V}\widehat{T^r1_{\{v\}}}(\xi)\,\omega^{x\cdot\xi}
 =\frac{1}{|V|}\sum_{\xi\in V}\omega^{r q(\xi)+(x-v)\cdot\xi}.
\]
Since \(q(\xi)=\xi\cdot\xi\), we have in \(\mathbb F_3\)
\[
 r q(\xi)+(x-v)\cdot\xi
 =r q\!\left(\xi+\frac{x-v}{2r}\right)-\frac{q(x-v)}{4r}.
\]
Substituting \(\eta=\xi+(x-v)/(2r)\) in the sum and expanding \(q(x-v)\) gives
\begin{equation}\label{eq:quadratic-kernel}
 \begin{aligned}
 (T^r1_{\{v\}})(x)
 &=\frac{\omega^{-q(x-v)/(4r)}}{|V|}\sum_{\eta\in V}\omega^{r q(\eta)}\\
 &=c_r |V|^{-1/2}\omega^{-q(x)/(4r)}\omega^{-q(v)/(4r)}
                    \omega^{x\cdot v/(2r)} .
 \end{aligned}
\end{equation}
Here \(c_r=|V|^{-1/2}\sum_{\eta\in V}\omega^{rq(\eta)}\) has modulus one, since
\(\left|1+2\omega^r\right|=\sqrt3\) and the sum factors over the coordinates.

We use normalized counting measure on \(V\), so \(\mathbb P(B)=|B|/3^n\) for every \(B\subseteq V\).
Fix \(r\in\{-1,1\}\) and write \(\gamma_v(x)=\omega^{x\cdot v/(2r)}\)
for the last factor in \eqref{eq:quadratic-kernel}. Thus
\(\gamma_{e_j}(x)=\omega^{x_j/(2r)}\). For every \(k=(k_1,\ldots,k_n)\in V\),
\[
 \mathbb P\bigl(\gamma_{e_j}(x)=\omega^{k_j}\text{ for all }1\le j\le n\bigr)
 =\mathbb P(x=2rk)=3^{-n}.
\]
Thus the \(\gamma_{e_j}\) are independent, each uniformly distributed on
\(\{1,\omega,\omega^2\}\).

Let \(S=\sum_{v\in I}a_v\gamma_v\), initially with \(\sum_{v\in I}|a_v|^2=1\).
Since \(1+\omega+\omega^2=0\), each \(a_v\gamma_v\) has mean zero.
The summands of \(\Re S\) are independent, and
\(-|a_v|\le\Re(a_v\gamma_v)\le|a_v|\). The same bounds hold for the
imaginary parts. Applying Hoeffding's inequality \cite[Theorem~2]{Hoeffding}
to each of the two sums at \(u/\sqrt2\), we obtain, for \(u>0\),
\[
 \mathbb P(|S|>u)
 \le \mathbb P(|\Re S|>u/\sqrt2)+\mathbb P(|\Im S|>u/\sqrt2)
 \le 2e^{-u^2/4}+2e^{-u^2/4}=4e^{-u^2/4}.
\]
For \(1<t<\infty\) it follows that
\[
 \|S\|_{L^t(V)}^t
 =t\int_0^\infty u^{t-1}\mathbb P(|S|>u)\,du
 \le 4t\int_0^\infty u^{t-1}e^{-u^2/4}\,du=:C_t^t .
\]
By homogeneity, for arbitrary coefficients,
\begin{equation}\label{eq:cube-root-sums}
 \Big\|\sum_{v\in I} a_v\gamma_v\Big\|_{L^t(V)}
 \le C_t\Big(\sum_{v\in I}|a_v|^2\Big)^{1/2}.
\end{equation}
Taking \(t=p'=p/(p-1)\) in \eqref{eq:cube-root-sums}, we obtain for \(f\in L^p(V)\)
\[
 \left(\sum_{v\in I}|\widehat f(\gamma_v)|^2\right)^{1/2}
 =\sup_{\sum_{v\in I}|a_v|^2=1}
 \left|\frac1{|V|}\sum_{x\in V} f(x)
              \overline{\sum_{v\in I}a_v\gamma_v(x)}\right|
 \le C_{p'}\|f\|_{L^p(V)} .
\]
Let \(\lambda_r(x)=\omega^{-q(x)/(4r)}\). By \eqref{eq:quadratic-kernel},
\(T^r1_{\{v\}}=c_r|V|^{-1/2}\omega^{-q(v)/(4r)}\lambda_r\gamma_v\).
Since \(T^r\) is unitary on \(L^2(V)\), \(T^rPT^{-r}\) is the orthogonal
projection onto the span of the orthonormal functions \(\lambda_r\gamma_v\), \(v\in I\). Thus
\[
 (T^rPT^{-r}f)(x)=\lambda_r(x)\sum_{v\in I}
       \widehat{\bar\lambda_r f}(\gamma_v)\gamma_v(x).
\]
Applying \eqref{eq:cube-root-sums} with \(t=p\) and then the coefficient
estimate to \(\bar\lambda_r f\) gives
\[
 \|T^rPT^{-r}f\|_{L^p(V)}
 \le C_p\left(\sum_{v\in I}
         |\widehat{\bar\lambda_r f}(\gamma_v)|^2\right)^{1/2}
 \le C_p C_{p'}\|f\|_{L^p(V)} .
\]
Define \(P_3:L^p(V^3)\to L^p(V^3)\) by \(P_3f=1_{I^3}f\).
Let \(\Pi:L^p(G)\to L^p(G)\) act on the sections of \(F\) by
\((\Pi F)_s=U_sP_3U_s^{-1}F_s\) for \(s=0,1,2\).
Since \(P_3=P\otimes P\otimes P\), the definition of \(U_s\) gives
\[
 U_sP_3U_s^{-1}
 = (T^{r_{s1}}PT^{-r_{s1}})\otimes(T^{r_{s2}}PT^{-r_{s2}})
   \otimes(T^{r_{s3}}PT^{-r_{s3}}).
\]
For each \(s\), the factor with \(r_{sj}=0\) is \(P\), of norm one. Each
of the other two factors has norm at most \(C_pC_{p'}\). Applying these
bounds in the corresponding coordinates and summing over \(s\) gives
\[
 \|\Pi F\|_{L^p(G)}^p
 =\frac13\sum_{s=0}^2\|U_sP_3U_s^{-1}F_s\|_{L^p(V^3)}^p
 \le (C_pC_{p'})^{2p}\|F\|_{L^p(G)}^p .
\]
Also \(\Pi^2=\Pi\) since \(P_3^2=P_3\). For \(F\in E\), \eqref{eq:quadratic-fibers} gives
\(\sum_s U_s^{-1}(\Pi F)_s=P_3\sum_s U_s^{-1}F_s=0\), so \(\Pi F\in E\).

Each function \(a\) in the range of \(P_3\) is determined by the \(n^3\) values
\(a_{ijk}=a(e_i,e_j,e_k)\), since
\[
 a(x)=
 \begin{cases}
 a_{ijk},&x=(e_i,e_j,e_k),\quad 1\le i,j,k\le n,\\
 0,&x\in V^3\setminus I^3 .
 \end{cases}
\]
With each index running from \(1\) to \(n\), define
\begin{equation}\label{eq:mixed-array-norms}
 \begin{aligned}
 &N_0(a)^p=\sum_i\Big(\sum_{j,k}|a_{ijk}|^2\Big)^{p/2},\\
 &N_1(a)^p=\sum_j\Big(\sum_{i,k}|a_{ijk}|^2\Big)^{p/2},\\
 &N_2(a)^p=\sum_k\Big(\sum_{i,j}|a_{ijk}|^2\Big)^{p/2}.
 \end{aligned}
\end{equation}
We shall prove that for every \(a=P_3a\) and \(s=0,1,2\),
\begin{equation}\label{eq:mixed-multiplier-norm}
 \|U_sa\|_{L^p(V^3)}\asymp_p |V|^{-1-1/p}N_s(a).
\end{equation}
The identity \(U_0=\mathrm{Id}\otimes T\otimes T^{-1}\) gives, for
\(1\le i\le n\) and \(x,y\in V\),
\[
 (U_0a)(e_i,x,y)
 =\sum_{j,k}a_{ijk}\,(T1_{\{e_j\}})(x)\,(T^{-1}1_{\{e_k\}})(y).
\]
Since \(q(e_j)=q(e_k)=1\), the two factors in \eqref{eq:quadratic-kernel}
that depend only on the input points \(e_j,e_k\) cancel. Taking moduli gives
\[
 |(U_0a)(e_i,x,y)|=|V|^{-1}|h_i(x,y)|,\qquad
 h_i(x,y)=\sum_{j,k=1}^n a_{ijk}\,\omega^{x_j/2}\omega^{-y_k/2}
 \quad(x,y\in V).
\]
The characters in this sum are orthonormal in \(L^2(V^2)\), so
\(\|h_i\|_{L^2(V^2)}^2=\sum_{j,k}|a_{ijk}|^2\).
Put \(b_{ij}(y)=\sum_k a_{ijk}\omega^{-y_k/2}\), so that
\(h_i(x,y)=\sum_j b_{ij}(y)\omega^{x_j/2}\).
For \(2\le t<\infty\), applying \eqref{eq:cube-root-sums} first in \(x\) and
then in \(y\), with the triangle inequality in \(L^{t/2}(V)\), gives
\[
 \begin{aligned}
 \|h_i\|_{L^t(V^2)}
 &\le C_t\left\|\sum_j |b_{ij}|^2\right\|_{L^{t/2}(V)}^{1/2}\\
 &\le C_t\left(\sum_j\|b_{ij}\|_{L^t(V)}^2\right)^{1/2}
 \le C_t^2\|h_i\|_{L^2(V^2)} .
 \end{aligned}
\]
Taking \(t=p\) proves the upper bound for \(p\ge2\), and the lower bound
follows from \(\|h_i\|_{L^2(V^2)}\le\|h_i\|_{L^p(V^2)}\).
When \(1<p<2\), we have \(\|h_i\|_{L^p(V^2)}\le\|h_i\|_{L^2(V^2)}\).
The reverse bound follows from H\"older's inequality and the estimate with \(t=p'>2\), since
\[
 \|h_i\|_{L^2(V^2)}^2
 \le \|h_i\|_{L^p(V^2)}\|h_i\|_{L^{p'}(V^2)}
 \le C_{p'}^2\|h_i\|_{L^p(V^2)}\|h_i\|_{L^2(V^2)} .
\]
Since \(U_0a\) vanishes when its first coordinate lies outside \(I\),
\[
 \|U_0a\|_{L^p(V^3)}^p
 =\frac{1}{|V|^{p+1}}\sum_i\|h_i\|_{L^p(V^2)}^p
 \asymp_p \frac{N_0(a)^p}{|V|^{p+1}} .
\]
For \(U_1=T^{-1}\otimes\mathrm{Id}\otimes T\) the outer sum is over \(j\),
giving \(N_1(a)^p\). For \(U_2=T\otimes T^{-1}\otimes\mathrm{Id}\)
it is over \(k\), giving \(N_2(a)^p\). This proves \eqref{eq:mixed-multiplier-norm}.

For \(F\in\Pi E\), put \(a_s=U_s^{-1}F_s\). Then \(P_3a_s=a_s\) for each \(s\) and \(a_0+a_1+a_2=0\) by \eqref{eq:quadratic-fibers}. Conversely these two conditions give \(\Pi F=F\) and \(F\in E\). Hence
\[
 \Pi E=\{(U_0a,U_1b,-U_2(a+b)):\operatorname{supp}a,\operatorname{supp}b\subset I^3\}.
\]
Identifying \(a,b\) with their arrays of values on \(I^3\) gives the graph space
\[
 Z=\{(a,b,a+b):a,b\in\mathbb C^{n\times n\times n}\},\qquad
 \|(a,b,a+b)\|_Z=N_0(a)+N_1(b)+N_2(a+b).
\]
For \(F=(U_0a,U_1b,-U_2(a+b))\), normalized counting measure in \(s\) gives
\[
 3\|F\|_{L^p(G)}^p=\|U_0a\|_{L^p(V^3)}^p+\|U_1b\|_{L^p(V^3)}^p+\|U_2(a+b)\|_{L^p(V^3)}^p .
\]
By \eqref{eq:mixed-multiplier-norm},
\[
 |V|^{1+1/p}\|F\|_{L^p(G)}
 \asymp_p \bigl(N_0(a)^p+N_1(b)^p+N_2(a+b)^p\bigr)^{1/p}.
\]
The \(\ell^p\) and \(\ell^1\) norms of these three numbers differ by a factor
of at most \(3^{1-1/p}\), so
\begin{equation}\label{eq:graph-isomorphism}
 |V|^{1+1/p}\|F\|_{L^p(G)}\asymp_p\|(a,b,a+b)\|_Z .
\end{equation}

\subsection{Lower bounds on the real coefficient spaces}
For any complex Banach space \(X\), write \(X_{\mathbb R}\) for the same
space with its given norm and scalar multiplication restricted to \(\mathbb R\).
Thus the triples in \(Z_{\mathbb R}\) still have complex entries.
We first show that \(\lust(Z_{\mathbb R})\le\lust(Z)\).
Fix a finite-dimensional real subspace \(F\subset Z_{\mathbb R}\). Its complex
span in \(Z\) is \(F+iF=\{f+ig:f,g\in F\}\). For any admissible constant
\(C\) for \(Z\), Definition~\ref{def:lust} applied to \(F+iF\) gives a complex
Banach space \(W\) with a \(1\)-unconditional basis \((w_j)\) and bounded
complex-linear operators
\[
 u:F+iF\longrightarrow W,\qquad v:W\longrightarrow Z,\qquad
 v(u(x))=x\quad(x\in F+iF),\qquad \|u\|\,\|v\|\le C.
\]
In \(W_{\mathbb R}\), the sequence \(w_1,iw_1,w_2,iw_2,\ldots\) is a real
\(1\)-unconditional basis, since changing the signs of the real and imaginary
parts of a coefficient preserves its modulus. Restricting \(u\) to \(F\) and
regarding \(v\) as real linear gives
\[
 F\xrightarrow{\ u|_F\ }W_{\mathbb R}\xrightarrow{\ v\ }Z_{\mathbb R},
 \qquad v(u(f))=f\quad(f\in F),\qquad \|u|_F\|\,\|v\|\le C.
\]
Thus every admissible constant for \(Z\) is also admissible for \(Z_{\mathbb R}\),
which proves the inequality. The same argument gives
\(\lust(X_{\mathbb R})\le\lust(X)\) for any complex Banach space \(X\).

Let \(Z_{\mathrm{re}}\) be the subspace of \(Z_{\mathbb R}\) consisting of
triples with real entries. The projection from \(Z_{\mathbb R}\) onto
\(Z_{\mathrm{re}}\) sends \((a,b,a+b)\) to \((\Re a,\Re b,\Re(a+b))\)
and has norm one by \eqref{eq:mixed-array-norms}. By \eqref{eq:graph-isomorphism},
\(Z\) is isomorphic to \(\Pi E\subset L^p(G)\) with distortion bounded in
terms of \(p\). Hence \eqref{eq:cotype} holds on its real subspaces with
exponent \(\max(2,p)\) and constant depending only on \(p\).

\noindent\textit{The case \(1<p<2\).}
Take \(Z_0=Z_{\mathrm{re}}\). Let \(e_{ijk}\) be the array whose only nonzero
entry is a one at \((i,j,k)\), and set
\[
 x_{ijk}=(0,e_{ijk},e_{ijk}),\qquad y_{ijk}=(e_{ijk},0,e_{ijk}),\qquad
 B_{ijk}=\operatorname{span}_{\mathbb R}\{x_{ijk},y_{ijk}\}.
\]
The coordinate expansion
\[
 (u,v,u+v)=\sum_{i,j,k=1}^n\bigl(v_{ijk}x_{ijk}+u_{ijk}y_{ijk}\bigr)
\]
is unique, so \(Z_0=\bigoplus_{i,j,k=1}^n B_{ijk}\). Multiplying both
\(u_{ijk}\) and \(v_{ijk}\) by \(\varepsilon_{ijk}\in\{-1,1\}\) preserves
each of \(N_0(u),N_1(v),N_2(u+v)\) in \eqref{eq:mixed-array-norms}.
Thus the decomposition into these two-dimensional blocks has unconditional
constant one. Since \(Z_0\) is finite dimensional and has cotype \(2\) with
constant depending only on \(p\), \cite[Proposition~2.2]{KT} gives an operator
\(L:Z_0\to Z_0\) preserving each \(B_{ijk}\) and satisfying
\begin{equation}\label{eq:block-operator}
 \|L\|\lesssim_p\lust(Z_0),\qquad
 \inf_{\theta\in\mathbb R}\|(L-\theta I)|_{B_{ijk}}\|\ge1/8 .
\end{equation}

Following the proof of \cite[Theorem~3.1]{KT}, write the matrix of
\(L|_{B_{ijk}}\) in the ordered basis \((x_{ijk},y_{ijk})\) as
\[
 \begin{pmatrix}a_{ijk}&b_{ijk}\\ c_{ijk}&d_{ijk}\end{pmatrix},
 \qquad \gamma_{ijk}=a_{ijk}-b_{ijk}+c_{ijk}-d_{ijk}.
\]
On this block \(\|s x_{ijk}+t y_{ijk}\|_Z=|s|+|t|+|s+t|\), which lies
between \(|s|+|t|\) and \(2(|s|+|t|)\). Hence
\[
 \begin{aligned}
 \|(L-d_{ijk}I)|_{B_{ijk}}\|
 &\le 2\max\{|a_{ijk}-d_{ijk}|+|c_{ijk}|,|b_{ijk}|\}\\
 &\le 8\max\{|b_{ijk}|,|c_{ijk}|,|\gamma_{ijk}|\}.
 \end{aligned}
\]
Taking \(\theta=d_{ijk}\) in \eqref{eq:block-operator} yields
\begin{equation}\label{eq:block-coefficients}
 \max\{|b_{ijk}|,|c_{ijk}|,|\gamma_{ijk}|\}\ge c_0,\qquad c_0=1/64 .
\end{equation}
By \eqref{eq:mixed-array-norms},
with the indices outside each sum fixed,
\[
 \Big\|\sum_i x_{ijk}\Big\|_Z
 =\Big\|\sum_j y_{ijk}\Big\|_Z
 =\Big\|\sum_k(x_{ijk}-y_{ijk})\Big\|_Z=2\sqrt n .
\]
Let \(\pi_1,\pi_2,\pi_3:Z_{\mathrm{re}}\to\mathbb R^{n\times n\times n}\)
be the three coordinate maps, so \(N_{s-1}(\pi_s z)\le\|z\|_Z\).
The block formulas for \(L\) give
\[
 \pi_1 Lx_{ijk}=c_{ijk}e_{ijk},\qquad
 \pi_2 Ly_{ijk}=b_{ijk}e_{ijk},\qquad
 \pi_3 L(x_{ijk}-y_{ijk})=\gamma_{ijk}e_{ijk}.
\]
For example,
\[
 \left(\sum_i|c_{ijk}|^p\right)^{1/p}
 =N_0\!\left(\pi_1 L\sum_i x_{ijk}\right)
 \le \|L\|\,\Big\|\sum_i x_{ijk}\Big\|_Z=2\sqrt n\,\|L\|.
\]
The same estimate for \(\pi_2 L\sum_j y_{ijk}\) and
\(\pi_3 L\sum_k(x_{ijk}-y_{ijk})\) gives
\[
 \left(\sum_j|b_{ijk}|^p\right)^{1/p}\le2\sqrt n\,\|L\|,
 \qquad
 \left(\sum_k|\gamma_{ijk}|^p\right)^{1/p}\le2\sqrt n\,\|L\|.
\]
Summing each estimate to the \(p\)th power over its two free indices and using \eqref{eq:block-coefficients} yields
\[
 c_0^p\le n^{-3}\sum_{ijk}(|b_{ijk}|^p+|c_{ijk}|^p+|\gamma_{ijk}|^p)
 \le 3\cdot 2^p n^{p/2-1}\|L\|^p .
\]
By \eqref{eq:block-operator}, it follows that \(\lust(Z_0)\gtrsim_p n^{1/p-1/2}\).

\noindent\textit{The case \(2<p<\infty\).}
Take \(n=m^2\) with \(m\) a positive integer.
Define \(Q:\mathbb R^{n\times n\times n}\to\mathbb R^{n\times n\times n}\) by
\[
 (Qa)_{ijk}=1_{\{k\le m\}}\frac1m\sum_{\ell=1}^m a_{ij\ell}.
\]
At fixed \(i,j\) this is an orthogonal projection in \(\ell^2\) over \(k\),
so \(N_0(Qa)\le N_0(a)\) and \(N_1(Qa)\le N_1(a)\).
For the third norm, let \(a_{\cdot\cdot k}=(a_{ijk})_{i,j}\). Convexity gives
\[
 N_2(Qa)^p
 =m\left\|\frac1m\sum_{k=1}^m a_{\cdot\cdot k}\right\|_2^p
 \le \sum_{k=1}^m\|a_{\cdot\cdot k}\|_2^p\le N_2(a)^p .
\]
Since \(Q^2=Q\), the map \((a,b,a+b)\mapsto(Qa,Qb,Qa+Qb)\)
is a norm-one projection on \(Z_{\mathrm{re}}\). Take \(Z_0\) to be its range.
Given \((a,b,a+b)\in Z_0\), define \(A,C\in\mathbb R^{n\times n}\) by
\(A_{ij}=\sqrt m\,a_{ij1}\) and \(C_{ij}=\sqrt m\,b_{ij1}\).
Since \(Qa=a\) and \(Qb=b\), we recover the arrays from these matrices by
\[
 a_{ijk}=m^{-1/2}A_{ij}1_{\{k\le m\}},\qquad
 b_{ijk}=m^{-1/2}C_{ij}1_{\{k\le m\}}.
\]
Conversely, these formulas give an element \(z(A,C)=(a,b,a+b)\in Z_0\)
for each pair of real \(n\times n\) matrices \(A,C\).
By \eqref{eq:mixed-array-norms}, with \(\rho=m^{1/p-1/2}\),
\[
 \begin{aligned}
 \|z(A,C)\|_Z
 &=\left(\sum_i\Big(\sum_j|A_{ij}|^2\Big)^{p/2}\right)^{1/p}
 +\left(\sum_j\Big(\sum_i|C_{ij}|^2\Big)^{p/2}\right)^{1/p}\\
 &\quad+\rho\left(\sum_{ij}|A_{ij}+C_{ij}|^2\right)^{1/2}.
 \end{aligned}
\]
Let \(e_{ij}\) be the matrix whose only nonzero entry is a one at \((i,j)\). Set
\[
 x_{ij}=z(0,e_{ij}),\qquad y_{ij}=z(e_{ij},0),\qquad
 B_{ij}=\operatorname{span}_{\mathbb R}\{x_{ij},y_{ij}\}.
\]
Each \(z(A,C)\) has the unique expansion
\[
 z(A,C)=\sum_{i,j=1}^n(C_{ij}x_{ij}+A_{ij}y_{ij}),
 \qquad Z_0=\bigoplus_{i,j=1}^n B_{ij}.
\]
Changing the signs of \(A_{ij}\) and \(C_{ij}\) together at any \((i,j)\)
preserves \(\|z(A,C)\|_Z\), so the decomposition has unconditional constant one.
Since \(Z_0\) has cotype \(p\) with constant depending only on \(p\),
\cite[Proposition~2.2]{KT} gives an operator \(L:Z_0\to Z_0\) preserving each
\(B_{ij}\) and satisfying \eqref{eq:block-operator} with \(B_{ij}\) in place of \(B_{ijk}\).
Write \(\left(\begin{smallmatrix}a_{ij}&b_{ij}\\c_{ij}&d_{ij}\end{smallmatrix}\right)\)
for the matrix of \(L|_{B_{ij}}\) in the ordered basis \((x_{ij},y_{ij})\),
and set \(\gamma_{ij}=a_{ij}-b_{ij}+c_{ij}-d_{ij}\). Since \(0<\rho\le1\),
\[
 |s|+|t|\le \|s x_{ij}+t y_{ij}\|_Z
 = |s|+|t|+\rho|s+t|\le 2(|s|+|t|).
\]
The argument for \eqref{eq:block-coefficients} therefore gives
\(\max\{|b_{ij}|,|c_{ij}|,|\gamma_{ij}|\}\ge c_0\) on each \(B_{ij}\).

Since \(n=m^2\) and \(\rho^2=n^{1/p-1/2}\), we have
\[
 \Big\|\sum_{\ell}x_{i\ell}\Big\|
 =\Big\|\sum_{\ell}y_{\ell j}\Big\|=n^{1/p}+\rho\sqrt n\le2\rho\sqrt n,
 \qquad
 \Big\|\sum_{i,j}(x_{ij}-y_{ij})\Big\|=2n\rho^2 .
\]
Substituting \(x_{ij}=z(0,e_{ij})\) and \(y_{ij}=z(e_{ij},0)\) into the block matrix formulas gives
\[
 \begin{aligned}
 &Lx_{ij}=a_{ij}x_{ij}+c_{ij}y_{ij}=z(c_{ij}e_{ij},a_{ij}e_{ij}),\\
 &Ly_{ij}=b_{ij}x_{ij}+d_{ij}y_{ij}=z(d_{ij}e_{ij},b_{ij}e_{ij}).
 \end{aligned}
\]
For their difference \(L(x_{ij}-y_{ij})=z(A,C)\), this gives
\(A+C=((c_{ij}-d_{ij})+(a_{ij}-b_{ij}))e_{ij}=\gamma_{ij}e_{ij}\).
Applying \(L\) to the three sums above and using the first, second and third
terms of \(\|z(A,C)\|_Z\), respectively, gives
\[
 \sum_j|c_{ij}|^2\le4\rho^2 n\|L\|^2,\quad
 \sum_i|b_{ij}|^2\le4\rho^2 n\|L\|^2,\quad
 \rho^2\sum_{ij}|\gamma_{ij}|^2\le4n^2\rho^4\|L\|^2 .
\]
Summing the first estimate over \(i\) and the second over \(j\), and using
the lower bound \(c_0\) on each block \(B_{ij}\), gives
\[
 c_0^2\le n^{-2}\sum_{ij}(|b_{ij}|^2+|c_{ij}|^2+|\gamma_{ij}|^2)
 \le12\rho^2\|L\|^2 .
\]
Thus \(\|L\|\gtrsim\rho^{-1}\), and the bound on \(\|L\|\) gives
\(\lust(Z_0)\gtrsim_p n^{(1/2-1/p)/2}\) for \(n=m^2\).

In either case \(Z_0\) is the range of a norm-one projection on \(Z_{\mathrm{re}}\),
and \(Z_{\mathrm{re}}\) is the range of a norm-one projection on \(Z_{\mathbb R}\). Applying \eqref{eq:transfer} to these
projections and then using the real factorization above gives
\(\lust(Z_0)\le\lust(Z_{\mathbb R})\le\lust(Z)\).
Recall that \(E=L_A^p(G)\). Applying \eqref{eq:graph-isomorphism},
the bound for \(\Pi|_E\), and \eqref{eq:spectrum-lust} in turn gives
\[
 \lust(Z)\lesssim_p\lust(\Pi E)\le\|\Pi|_E\|\,\lust(E)
 \lesssim_p\lust\bigl(H^p(\mathbb T^\infty)\bigr).
\]
Letting \(n\to\infty\), through squares when \(p>2\), gives
\(\lust(H^p(\mathbb T^\infty))=\infty\) for \(1<p<\infty\), \(p\ne2\).

\subsection{Local unconditional structure at the endpoint}\label{sec:quotient}

For a real Banach space \(Y\) put
\[
 d_n(Y)=\inf\{\|i\|\|i^{-1}\|:\ i:\ell_\infty^n(\mathbb R)\to Y
                  \text{ is a linear embedding}\},
\]
where the inverse is defined on \(i(\ell_\infty^n)\) and \(\inf\varnothing=\infty\).
Johnson calls the property in Definition~\ref{def:lust} GL-l.u.st.
By \cite[Corollary~2.2]{Johnson}, if \(q:L\to Y\) is a bounded
linear surjection and \(L\) embeds isomorphically in a Banach lattice \(B\)
with \(\sup_n d_n(B)=\infty\), then
\begin{equation}\label{eq:johnson-kernel}
 \ker q\ \text{has local unconditional structure}\quad\Longrightarrow\quad
                    \sup_n d_n(Y)=\infty .
\end{equation}

For each \(n\ge1\), take \(G=(\mathbb Z/4\mathbb Z)^n\).
Put \(\chi_j(x)=i^{x_j}\) for \(1\le j\le n\) and
\(A=\widehat G\setminus\{\chi_1,\ldots,\chi_n\}\).
Write \([f]=f+L_A^1(G)\), with norm \(\inf_{h\in L_A^1(G)}\|f-h\|_{L^1(G)}\).
By quotient duality, \(\|[f]\|\) is the supremum of
\(\bigl|\frac1{|G|}\sum_{x\in G}f(x)\overline{g(x)}\bigr|\) over
\(\|g\|_{L^\infty(G)}\le1\) in the annihilator
\[
 (L_A^1(G))^\perp
 =\{g\in L^\infty(G):\widehat g(\chi)=0\ \text{for all }\chi\in A\}
 =\Big\{\sum_{j=1}^n a_j\chi_j:a_j\in\mathbb C\Big\}.
\]
For each \(j\), choose \(x_j\in\mathbb Z/4\mathbb Z\) to maximize \(\Re(a_j i^{x_j})\).
Then \(\Re(a_j i^{x_j})=\max\{|\Re a_j|,|\Im a_j|\}\ge |a_j|/\sqrt2\). Consequently
\[
 2^{-1/2}\sum_j|a_j|\le\left\|\sum_ja_j\chi_j\right\|_\infty\le\sum_j|a_j|.
\]
Substituting \(g=\sum_j a_j\chi_j\) into the pairing therefore gives
\[
 \max_j|\widehat f(\chi_j)|\le\|[f]\|\le\sqrt2\,\max_j|\widehat f(\chi_j)|.
\]
The map \([f]\mapsto(\widehat f(\chi_j))_{j=1}^n\) is onto
\(\ell_\infty^n(\mathbb C)\), since \(f=\sum_j b_j\chi_j\) gives the
vector \((b_j)_{j=1}^n\).

Apply Theorem~\ref{thm:retract} with \(p=1,\ \varepsilon=1\). Since
\(J(L_A^1(G))\subseteq H^1(\mathbb T^\infty)\) and
\(R(H^1(\mathbb T^\infty))\subseteq L_A^1(G)\), \(J\) and \(R\) induce
\[
\begin{aligned}
 &\bar J:L^1(G)/L_A^1(G)\longrightarrow L^1(\mathbb T^\infty)/H^1(\mathbb T^\infty),
 &&\bar J[f]=[Jf],\\
 &\bar R:L^1(\mathbb T^\infty)/H^1(\mathbb T^\infty)\longrightarrow L^1(G)/L_A^1(G),
 &&\bar R[F]=[RF].
 \end{aligned}
\]
Here \([F]=F+H^1(\mathbb T^\infty)\) for \(F\in L^1(\mathbb T^\infty)\).
Taking infima over representatives gives \(\|\bar J\|\le\|J\|\) and \(\|\bar R\|\le\|R\|\).
Since \(\bar R\bar J[f]=[RJf]=[f]\), the inverse of \(\bar J\) on its
range is the restriction of \(\bar R\). Thus \eqref{eq:retraction-bounds}
gives \(\|\bar J\|\,\|\bar J^{-1}\|\le2\). Composing \(\bar J\) with the
inverse of the coefficient map above embeds \(\ell_\infty^n(\mathbb C)\)
in \(L^1(\mathbb T^\infty)/H^1(\mathbb T^\infty)\) with distortion at most \(2\sqrt2\).

To apply \eqref{eq:johnson-kernel}, take \(L=L^1(\mathbb T^\infty)_{\mathbb R}\)
and \(B=L^1(\mathbb T^\infty;\mathbb R)\oplus_1 L^1(\mathbb T^\infty;\mathbb R)\).
The real-linear map \(f\mapsto(\Re f,\Im f)\) embeds \(L\) into \(B\), since
\(\|f\|_1\le\|\Re f\|_1+\|\Im f\|_1\le\sqrt2\,\|f\|_1\).
The lattice \(B\) is the real \(L^1\) space of two disjoint copies of the
torus and has cotype two by \cite[\S1.e, p.~73]{LT}.
For any embedding \(i:\ell_\infty^n(\mathbb R)\to B\), let \(e_1,\ldots,e_n\)
be the standard unit vectors of \(\ell_\infty^n(\mathbb R)\).
With \(C_2(B)\) a cotype-two constant for \(B\) in \eqref{eq:cotype}, we obtain
\[
 \frac{\sqrt n}{\|i^{-1}\|}
 \le \left(\sum_{j=1}^n\|i(e_j)\|^2\right)^{1/2}
 \le C_2(B)\left(\mathbb E\left\|i\left(\sum_{j=1}^n\varepsilon_j e_j\right)\right\|^2\right)^{1/2}
 \le C_2(B)\|i\|,
\]
since \(\|\sum_{j=1}^n\varepsilon_j e_j\|_\infty=1\) for every choice of signs.
Thus \(d_n(B)\ge\sqrt n/C_2(B)\). Let
\(Y=(L^1(\mathbb T^\infty)/H^1(\mathbb T^\infty))_{\mathbb R}\).
Restricting the preceding embeddings to \(\ell_\infty^n(\mathbb R)\) gives
\(d_n(Y)\le2\sqrt2\) for every \(n\). The quotient map \(q:L\to Y\),
\(q(f)=[f]\), has kernel \(H^1(\mathbb T^\infty)_{\mathbb R}\).
Since \(\sup_n d_n(B)=\infty\) and \(\sup_n d_n(Y)<\infty\),
\eqref{eq:johnson-kernel} gives \(\lust(H^1(\mathbb T^\infty)_{\mathbb R})=\infty\).
The inequality \(\lust(X_{\mathbb R})\le\lust(X)\) proved above yields
\(\lust(H^1(\mathbb T^\infty))=\infty\).

The Bohr isometry \(\mathcal H^p\simeq H^p(\mathbb T^\infty)\) now gives
\(\lust(\mathcal H^p)=\infty\) for \(p\ne2\), so Lemma~\ref{lem:lust-transfer}
excludes an unconditional basis in this range. For \(p=2\), the Dirichlet
monomials form an orthonormal basis and \(\lust(\mathcal H^2)=1\) by the
same lemma. This proves Theorem~\ref{thm:unc}.\qed

\section{Proof of Theorem~\ref*{thm:basis}}\label{sec:basis-proof}

Put \(X=H^1(\mathbb T^\infty)\). We first construct a Schauder basis of \(X\) by enlarging its finite-dimensional blocks. Section~\ref{sec:explicit} gives a choice of blocks and maps for which both the basis vectors and their coefficient functionals have finite formulas in the Dirichlet series variables.

\subsection{An unconditional decomposition with blocks isomorphic to polynomial spaces}
Write \(q(\alpha)=\sum_{j\ge1}j\alpha_j\) for \(\alpha\in\mathbb N_0^{(\mathbb N)}\).
For \(m\ge1\) let \(F_m=\operatorname{span}\{z^\alpha:q(\alpha)\le m\}\subset X\).
This space has finite dimension, since \(q(\alpha)\le m\) implies
\(\alpha_j=0\) for \(j>m\) and \(\alpha_j\le m/j\). For a function
\(a:[0,\infty)\to\mathbb C\), define \(a(q)\) on analytic polynomials by
\(a(q)z^\alpha=a(q(\alpha))z^\alpha\), extended linearly.

Take a smooth function \(\eta\) on \([0,\infty)\) equal to one on \([0,1]\) and to zero on \([2,\infty)\), and set
\[
 \phi_0(x)=\eta(x),\qquad
 \phi_k(x)=\eta(2^{-k}x)-\eta(2^{1-k}x)\quad(k\ge1),\qquad m_k=2^{k+1}\quad(k\ge0).
\]
Then \(\sum_{k=0}^N\phi_k(x)=\eta(2^{-N}x)\to1\). Let \(\psi_k(x)=\psi(2^{-k}x)\) for \(k\ge1\), where \(\psi\) is smooth, supported in \([1/4,4]\) and equal to one on \([1/2,2]\). Take \(\psi_0\) smooth, compactly supported and equal to one on \([0,2]\). Thus \(\psi_k\phi_k=\phi_k\).
For \(a_k=\phi_k\) or \(a_k=\psi_k\), the support of \(a_k\) for \(k\ge1\)
is contained in \([2^{k-2},2^{k+2}]\). Each \(x>0\) belongs to at most
five of these intervals, on each of which \(2^{-k}x\le4\). Hence every
finite sum \(h=\sum_k\varepsilon_ka_k\), \(\varepsilon_k\in\{0,\pm1\}\),
satisfies \(\sup_{x\ge0}|h(x)|+\sup_{x>0}|xh'(x)|\le C\),
where \(C\) depends only on the fixed functions \(\eta,\psi,\psi_0\).
Set \(m(t)=h(\log_2 t)\) for \(t\ge1\) and \(m(t)=h(0)\) for \(0<t<1\).
Then \(m(2^\ell)=h(\ell)\) for \(\ell\ge0\), and
\(t\log(t)m'(t)=(\log_2 t)h'(\log_2 t)\) for \(t>1\).
For every analytic polynomial \(\sum_\ell b_\ell w^\ell\), the Bohr isometry
and the multiplier bound of Aleman, Olsen and Saksman \cite{AAOS} give
\[
 \begin{aligned}
 \left\|\sum_\ell h(\ell)b_\ell w^\ell\right\|_{L^1(\mathbb T)}
 &=\left\|\sum_\ell m(2^\ell)b_\ell 2^{-\ell s}\right\|_{\mathcal H^1}\\
 &\lesssim
 \left(\|m\|_\infty+\sup_{t>e}|t\log(t)m'(t)|\right)
 \left\|\sum_\ell b_\ell 2^{-\ell s}\right\|_{\mathcal H^1}\\
 &\lesssim\left\|\sum_\ell b_\ell w^\ell\right\|_{L^1(\mathbb T)} .
 \end{aligned}
\]
For \(f(z)=\sum_\alpha c_\alpha z^\alpha\), apply this estimate in \(w\) to
\(f(wz_1,w^2z_2,\ldots)=\sum_\alpha c_\alpha z^\alpha w^{q(\alpha)}\).
Integrating in \(z\) gives
\begin{equation}\label{eq:block-multiplier}
 \left\|\sum_k\varepsilon_k a_k(q)f\right\|_X
 =\left\|\sum_\alpha h(q(\alpha))c_\alpha z^\alpha\right\|_X
 \le C\|f\|_X ,
\end{equation}
where \(C\) is independent of the finite choice of signs.

Let \(Y\) be the completion of finite arrays \(g=(g_{kl})_{k\ge0,l\ge1}\), \(g_{kl}\in F_{m_k}\), for the norm
\begin{equation}\label{eq:block-array-norm}
 \|g\|_Y=\int_{\mathbb T^\infty}\left(\sum_{k,l}|g_{kl}(z)|^2\right)^{1/2}dm(z).
\end{equation}
For an analytic polynomial \(f\in X\) and a finite array \(g\in Y\), define \(Af\in Y\) and \(Bg\in X\) by
\[
 (Af)_{kl}=1_{\{l=1\}}\phi_k(q)f,\qquad Bg=\sum_k\psi_k(q)g_{k1}.
\]
Since \(\phi_k(x)=0\) for \(x>m_k\), each \(\phi_k(q)f\) belongs to \(F_{m_k}\).
Let \((\varepsilon_k)\) be independent Rademacher random variables.
Applying \eqref{eq:unit-moment} at each \(z\) to the coefficients
\((\phi_k(q)f)(z)\), followed by \eqref{eq:block-multiplier} with \(a_k=\phi_k\), gives
\(\|Af\|_Y\lesssim\mathbb E\|\sum_k\varepsilon_k\phi_k(q)f\|_1\lesssim\|f\|_1\).
For a finite array \(g\), take both sums below over the indices with
\(g_{k1}\ne0\). Since \(\mathbb E(\varepsilon_k\varepsilon_j)=1_{\{k=j\}}\),
\[
 \mathbb E\left[\left(\sum_k\varepsilon_k\psi_k(q)\right)
                    \left(\sum_j\varepsilon_j g_{j1}\right)\right]
 =\sum_{k,j}\mathbb E(\varepsilon_k\varepsilon_j)\psi_k(q)g_{j1}
 =\sum_k\psi_k(q)g_{k1}=Bg .
\]
Using \eqref{eq:block-multiplier} with \(a_k=\psi_k\) and then \eqref{eq:unit-moment}, we obtain
\[
 \|Bg\|_1\lesssim\mathbb E\left\|\sum_j\varepsilon_jg_{j1}\right\|_1
 \lesssim \int_{\mathbb T^\infty}\left(\sum_j|g_{j1}(z)|^2\right)^{1/2}dm(z)
 \le\|g\|_Y .
\]
Thus \(A:X\to Y\) and \(B:Y\to X\) extend boundedly. On monomials,
\[
 BAz^\alpha=\sum_k\psi_k(q(\alpha))\phi_k(q(\alpha))z^\alpha=z^\alpha .
\]
Hence \(BA=I_X\) by density, and \(AB:Y\to Y\) is a projection onto \(A(X)\).

To embed \(Y\) in \(X\), use product variables
\((t,z,u)\in\mathbb T\times\mathbb T^\infty\times\mathbb T^\infty\),
with \(z=(z_j)_{j\ge1}\) and \(u=(u_{kl})_{k\ge0,l\ge1}\).
For a finite array \(g\) set
\begin{equation}\label{eq:array-embedding}
 (Jg)(t,z,u)=\sum_{kl}u_{kl}t^{m_k}g_{kl}(z_1/t,z_2/t^2,\ldots).
\end{equation}
If \(g\) has the single entry \(g_{kl}=z^\alpha\), then
\(Jg=u_{kl}t^{m_k-q(\alpha)}z^\alpha\). As \(q(\alpha)\le m_k\), this is an analytic monomial.
For fixed \(z,t\), \eqref{eq:unit-moment} applied to the independent circle
coordinates \(u_{kl}\) with product Haar measure gives
\[
 \int_{\mathbb T^\infty} |(Jg)(t,z,u)|\,dm(u)
 \asymp \left(\sum_{k,l}|g_{kl}(z_1/t,z_2/t^2,\ldots)|^2\right)^{1/2}.
\]
After integration in \(z,t\), invariance of \(dm(z)\) under
\(z_j\mapsto z_j/t^j\) yields \(\|Jg\|_1\asymp\|g\|_Y\). Thus \(J\) extends to an isomorphism
of \(Y\) onto a closed subspace of \(X\). Define two contractions \(P_1,P_2:X\to X\) by
\[
 \begin{aligned}
 &(P_1 F)(t,z,u)=\int_{\mathbb T}\bar\lambda F(t,z,\lambda u)\,dm(\lambda),\\
 &(P_2 F)(t,z,u)=\int_{\mathbb T}
 F(\lambda t,(\lambda^jz_j)_j,(\lambda^{-m_k}u_{kl})_{kl})\,dm(\lambda).
 \end{aligned}
\]
For an analytic monomial \(z^\alpha t^a u^\beta\), let \(|\beta|=\sum_{k,l}\beta_{kl}\). Using \(\int_{\mathbb T}\lambda^d\,dm(\lambda)=1_{\{d=0\}}\), we obtain
\[
 P_2P_1(z^\alpha t^a u^\beta)=
 \begin{cases}
 z^\alpha t^a u_{kl}, & \beta=e_{kl}\ \text{and }a=m_k-q(\alpha),\\
 0,&\text{otherwise},
 \end{cases}
\]
where \(e_{kl}\) is the exponent array with a single one at \((k,l)\).
Thus \((P_2P_1)^2=P_2P_1\) on analytic polynomials and, by density, on \(X\). By \eqref{eq:array-embedding}, its range is
\[
 \overline{\operatorname{span}}\{z^\alpha t^{m_k-q(\alpha)}u_{kl}:
                           k\ge0,\ l\ge1,\ q(\alpha)\le m_k\}=J(Y).
\]

For a Banach space \(Z\), let \(\operatorname{Rad}(Z)\) be the completion of the finitely supported sequences \((z_i)_{i\ge1}\), \(z_i\in Z\), under the norm \(\mathbb E\|\sum_i\varepsilon_i z_i\|_Z\), where the \(\varepsilon_i\) are independent Rademacher variables.
Independence of the signs gives, for each finite sum \(S=\sum_{i=1}^N\varepsilon_i z_i\),
\[
 \mathbb E(S\mid\varepsilon_1)=\varepsilon_1z_1,\qquad
 \mathbb E(S\mid\varepsilon_2,\ldots,\varepsilon_N)=\sum_{i=2}^N\varepsilon_i z_i .
\]
Both coordinate projections are therefore contractive. Splitting off the
first coordinate gives \(\operatorname{Rad}(Z)\simeq Z\oplus\operatorname{Rad}(Z)\),
where \(\simeq\) denotes Banach-space isomorphism and the external direct sum
is equipped with the sum norm.
If \(L:Z\to Z'\) is an isomorphism, applying \(L\) to each coordinate
gives \(\operatorname{Rad}(Z)\simeq\operatorname{Rad}(Z')\).
For \(Z=Z_1\oplus_1 Z_2\), the expectation of the sum norm is the sum of
the two expectations, so \(\operatorname{Rad}(Z)=
\operatorname{Rad}(Z_1)\oplus_1\operatorname{Rad}(Z_2)\) under the coordinate identification.
For a finite sequence of finite arrays \(g^{(i)}=(g^{(i)}_{kl})_{k,l}\in Y\),
apply \eqref{eq:unit-moment} to the vectors \((g^{(i)}_{kl}(z))_{k,l}\in\ell^2\) and integrate in \(z\) to obtain
\[
 \|(g^{(i)})_i\|_{\operatorname{Rad}(Y)}
 \asymp\int_{\mathbb T^\infty}\left(\sum_{i,k,l}|g^{(i)}_{kl}(z)|^2\right)^{1/2}dm(z).
\]
For fixed \(k\), all entries belong to the same \(F_{m_k}\), so reindexing \((i,l)\) by a bijection \(\mathbb N^2\to\mathbb N\) gives \(\operatorname{Rad}(Y)\simeq Y\) by \eqref{eq:block-array-norm}.
The projections \(P_2P_1:X\to J(Y)\) and \(AB:Y\to A(X)\) give
\(X\simeq Y\oplus W\) and \(Y\simeq X\oplus V\), respectively, with
\(W=\ker(P_2P_1)\) and \(V=\ker B\). Applying \(\operatorname{Rad}\)
to \(Y\simeq Y\oplus W\oplus V\) and using \(\operatorname{Rad}(Y)\simeq Y\) gives
\[
 Y\simeq Y\oplus\operatorname{Rad}(W)\oplus\operatorname{Rad}(V).
\]
Since \(W\oplus\operatorname{Rad}(W)\simeq\operatorname{Rad}(W)\), adjoining
\(W\) to the right-hand side leaves it isomorphic to \(Y\). Thus \(X\simeq Y\).
Fix an isomorphism \(H:Y\to X\) and enumerate the coordinate pairs \((k,l)\).
If \((k,l)\) is the \(n\)th pair, let \(D_n\) be the image under \(H\) of the arrays
supported at \((k,l)\). By \eqref{eq:block-array-norm}, changing the signs of
whole coordinates preserves the norm of \(Y\), and the finite arrays are dense.
Thus \((D_n)\) is an unconditional finite-dimensional decomposition of \(X\).
The map \(d\mapsto(H^{-1}d)_{kl}\) identifies \(D_n\) with \(F_{m_k}\) with
distortion at most \(\|H\|\|H^{-1}\|\), independently of \(n\).

\subsection{Enlarging the blocks and ordering their bases}\label{sec:blocks}

We first construct a complemented copy of \(\ell^1(X)\) in \(X\), then use
the direct-sum argument in the proof of Pe{\l}czy\'nski's Proposition~4
\cite[pp.~212--213]{Pelczynski} to absorb the spaces needed to enlarge the blocks.

Use product variables \((t,z)\in\mathbb T\times\mathbb T^\infty\) for \(X\), with \(z=(z_j)_{j\ge1}\). Take pairwise disjoint open arcs \(I_j\) with centers \(\zeta_j\in\mathbb T\). Fix \(0<\delta<1/2\). For each \(j\) take
\[
 h_j(t)=\frac{(1+\bar\zeta_jt)^{N_j}}
                  {\|(1+\bar\zeta_jt)^{N_j}\|_{L^1(\mathbb T)}}
 \quad\text{with}\quad
 \int_{\mathbb T\setminus I_j}|h_j(t)|\,dm(t)<\delta/2 .
\]
Such an \(N_j\) exists because \(|1+\bar\zeta_jt|\) is strictly larger
on a sufficiently small closed arc about \(\zeta_j\) than anywhere outside \(I_j\).
Put \(v_j=h_j1_{I_j}/\int_{I_j}|h_j|\,dm\). Then
\(\|h_j-v_j\|_{L^1(\mathbb T)}=2\int_{\mathbb T\setminus I_j}|h_j|\,dm<\delta\).
Define
\[
 \begin{aligned}
 &S:\ell^1(X)\longrightarrow X,
 &&[S(f_j)](t,z)=\sum_j h_j(t)f_j(z),\\
 &C:X\longrightarrow\ell^1(X),
 &&(CF)_j(z)=\int_{I_j}F(t,z)\,\overline{\operatorname{sgn}v_j(t)}\,dm(t).
 \end{aligned}
\]
The same formula extends \(C\) to a contraction from \(L^1(\mathbb T\times\mathbb T^\infty)\) to \(\ell^1(L^1(\mathbb T^\infty))\). As the arcs are disjoint,
\[
 \sum_j\|(CF)_j\|_{L^1(\mathbb T^\infty)}
 \le\sum_j\int_{I_j}\|F(t,\cdot)\|_{L^1(\mathbb T^\infty)}\,dm(t)
 \le\|F\|_{L^1(\mathbb T\times\mathbb T^\infty)}.
\]
For analytic polynomials each \((CF)_j\) is analytic in \(z\), so \(C(X)\subset\ell^1(X)\) by density.
As \(\int_{I_i}v_j\overline{\operatorname{sgn}v_i}\,dm=1_{\{i=j\}}\), we have
\[
 \|CS(f_j)-(f_j)\|_{\ell^1(X)}
 \le \left\|\sum_j (h_j-v_j)f_j\right\|_{L^1(\mathbb T\times\mathbb T^\infty)}
 \le \delta\sum_j\|f_j\|_X .
\]
Thus \(CS\) is invertible and \((CS)^{-1}C:X\to\ell^1(X)\) is a bounded left inverse of \(S\).

With \(B=\ell^1(X)\) and \(W=\ker C\), the projection \(S(CS)^{-1}C\) gives
\(X\simeq B\oplus W\). By \cite[Prop.~3(b), p.~212]{Pelczynski},
\[
 B=\ell^1(X)\simeq\ell^1(B)\oplus_1\ell^1(W)
 \simeq B\oplus_1\ell^1(W).
\]
Here \(\ell^1(B)\simeq B\) by reindexing \(\mathbb N^2\) by \(\mathbb N\).
Since \(W\oplus_1\ell^1(W)\simeq\ell^1(W)\), this gives \(W\oplus B\simeq B\)
and hence \(X\simeq\ell^1(X)\). More generally, if \(X\simeq V\oplus Z\), then
\begin{equation}\label{eq:absorb}
 X\simeq\ell^1(V)\oplus\ell^1(Z)\simeq V\oplus X ,
\end{equation}
since \(V\oplus_1\ell^1(V)\simeq\ell^1(V)\).

\paragraph{Bases for the enlarged blocks.}
The monomials spanning \(F_m\) are indexed by \(\Lambda_m=\{\alpha\in\mathbb N_0^m:\sum_{j=1}^m j\alpha_j\le m\}\). For integers \(L>m\) put \(G_L=(\mathbb Z/L\mathbb Z)^m\) and \(\zeta_L=e^{2\pi i/L}\). Define \(S_L:F_m\to L^1(G_L)\) by
\[
 (S_Lf)(x)=f(\zeta_L^{x_1},\ldots,\zeta_L^{x_m}),\qquad
 S_L(z^\alpha)(x)=\zeta_L^{\sum_j x_j\alpha_j}.
\]
Since \(0\le\alpha_j\le m<L\), distinct \(\alpha\in\Lambda_m\) give distinct characters. Thus \(S_L F_m=L^1_{A_L}(G_L)\), where \(A_L=\{S_L(z^\alpha):\alpha\in\Lambda_m\}\).
The Riemann sums \(L^{-m}\sum_{x\in G_L}|S_L f(x)|\) converge to \(\|f\|_{L^1(\mathbb T^m)}\).
On the \(L^1\) unit sphere of \(F_m\) this convergence is uniform, since the
sphere is compact in the uniform norm and the Riemann sums are \(1\)-Lipschitz in that norm.
Choose \(L=L(m)>m\) so that for every \(f\in F_m\),
\[
 \tfrac23\|f\|_{L^1(\mathbb T^m)}\le\|S_Lf\|_{L^1(G_L)}
 \le\tfrac43\|f\|_{L^1(\mathbb T^m)} .
\]

Write \((k_1,l_1),(k_2,l_2),\ldots\) for the enumeration of the array
coordinates used to define the blocks \(D_1,D_2,\ldots\).
Fix \(n\ge1\). For \(d\in D_n\), the array \(H^{-1}d\) is supported at
\((k_n,l_n)\). Set \(m=m_{k_n}\) and take \(L=L(m)\) as above. The coordinate isomorphism is
\[
 H_n:D_n\longrightarrow F_m,\qquad H_nd=(H^{-1}d)_{k_n l_n}.
\]
Take \(G=G_L\), \(A=A_L\), \(p=1\), and \(\varepsilon=1\) in the proof of
Theorem~\ref{thm:retract}. With \(r\) chosen in Step~1 of that proof,
the index set \(D\) from Step~3 for the block \(D_n\) is
\[
 D(n):=\left\{\beta\in\{0,1\}^{m\times r}:
 \left(\sum_{j=1}^r\beta_{ij}\right)_{i=1}^m\equiv\alpha\pmod L
 \text{ for some }\alpha\in\Lambda_m\right\}.
\]
Let \(E_{D(n)}=\operatorname{span}\{w^\beta:\beta\in D(n)\}\) with the
norm inherited from \(L^1(\mathbb T^{mr})\). The maps \(S,Q\) from
Steps~1 and~2 of that proof restrict to
\[
 S:L^1_{A_L}(G_L)\to E_{D(n)},\qquad
 Q:E_{D(n)}\to L^1_{A_L}(G_L),
 \qquad QS=I,\quad\|S\|\|Q\|\le2
\]
by \eqref{eq:finite-left-inverse} and the inclusions in that proof. Hence the maps
\[
 i_n=SS_LH_n:D_n\to E_{D(n)},\qquad
 r_n=H_n^{-1}S_L^{-1}Q:E_{D(n)}\to D_n
\]
satisfy \(r_ni_n=I_{D_n}\) and
\(\|i_n\|\|r_n\|\le4\|H\|\|H^{-1}\|\), where \(S_L^{-1}\) acts on \(L^1_{A_L}(G_L)\).
Replace \(i_n,r_n\) by \(i_n/\|i_n\|,\ \|i_n\|r_n\), respectively, and let \(V_n=\ker r_n\). Then
\[
 (d,v)\mapsto i_nd+v,\qquad
 e\mapsto(r_ne,e-i_nr_ne)
\]
are inverse isomorphisms between \(D_n\oplus_1V_n\) and \(E_{D(n)}\) with uniform bounds.
For this choice of \(D(n)\), take \(U,V\) from Step~4 and \(\iota,P\) from
Step~5 of the proof of Theorem~\ref{thm:retract}. Then \(\iota U:E_{D(n)}\to X\)
is an isometry and \(VP:X\to E_{D(n)}\) is a contraction with
\((VP)(\iota U)=I_{E_{D(n)}}\). Consequently
\[
 (\iota U)(I-i_nr_n)(VP):X\longrightarrow X
\]
is a projection onto \((\iota U)(V_n)\) of norm at most \(1+\|i_n\|\|r_n\|\).
Applying these projections coordinatewise in \(\ell^1(X)\), as in
\cite[Prop.~3(d), p.~212]{Pelczynski}, gives a complemented copy of
\(V=(\bigoplus_n V_n)_1\) in \(\ell^1(X)\simeq X\).

By \eqref{eq:absorb}, \(X\simeq X\oplus_1 V\).
For \(x=\sum_n d_n\) with \(d_n\in D_n\) and \(v=(v_n)_n\in V\), define
\[
 P_j(x,v)=\left(\sum_{n=1}^j d_n,(v_1,\ldots,v_j,0,\ldots)\right),\qquad P_0=0.
\]
Keeping the first \(j\) coordinates of \(H^{-1}x\) in the enumeration of
the pairs \((k,l)\) can only decrease the sum under the square root in \eqref{eq:block-array-norm}.
Applying \(H\) to this truncated array gives \(\sum_{n=1}^j d_n\). Hence
\[
 \|P_j(x,v)\|_{X\oplus_1V}
 =\Big\|\sum_{n=1}^j d_n\Big\|_X+\sum_{n=1}^j\|v_n\|_{V_n}
 \le \|H\|\|H^{-1}\|\|x\|_X+\|v\|_V .
\]
Also \(P_j(x,v)\to(x,v)\), since \(\sum_n d_n=x\) in \(X\) and
\(\sum_n\|v_n\|_{V_n}<\infty\).

In block \(j\), order \(D(j)\) by \eqref{eq:binary-order}, taking the
matrix entry \(\beta_{i\ell}\) as coordinate \((i-1)r+\ell\).
Write its elements in this order as \(\beta_1,\ldots,\beta_{|D(j)|}\).
Let \(T_j:E_{D(j)}\to D_j\oplus_1V_j\) denote the isomorphism
\(T_je=(r_je,e-i_jr_je)\). Since the monomials \(w^{\beta_\ell}\)
form a basis of \(E_{D(j)}\), their images \((T_jw^{\beta_\ell})_{\ell=1}^{|D(j)|}\)
form an ordered basis of block \(j\).
For \(e=\sum_{\ell=1}^{|D(j)|}a_\ell w^{\beta_\ell}\), the projection
\(B_{j,k}\) onto the first \(k\) vectors of this basis satisfies
\[
 B_{j,k}(T_je)=T_j\Big(\sum_{\ell=1}^k a_\ell w^{\beta_\ell}\Big),
 \qquad
 \|B_{j,k}(T_je)\|_{D_j\oplus_1V_j}
 \le C\|T_j\|\|T_j^{-1}\|\|T_je\|_{D_j\oplus_1V_j}.
\]
Here \(C\) is the absolute constant in \eqref{eq:binary-prefix-bound}.
The uniform bounds on \(T_j\) and \(T_j^{-1}\) therefore give
\(\sup_{j,k}\|B_{j,k}\|<\infty\).

List the block bases consecutively and let \(S_N\) denote the partial sum
through the \(N\)th vector. There are unique \(j\ge1\) and \(1\le k\le |D(j)|\)
such that \(N=\sum_{n<j}|D(n)|+k\). The first \(j-1\) blocks contribute
\(P_{j-1}(x,v)\). The \(j\)th block is \((P_j-P_{j-1})(x,v)\), identified
with \((d_j,v_j)\in D_j\oplus_1V_j\), and contributes its first \(k\) terms. Thus
\begin{equation}\label{eq:concatenated-partial-sum}
 S_N=P_{j-1}+B_{j,k}(P_j-P_{j-1}).
\end{equation}
In particular, \(\sup_N\|S_N\|\le
(\sup_j\|P_j\|)(1+2\sup_{j,k}\|B_{j,k}\|)<\infty\). For each \((x,v)\in X\oplus_1V\),
\[
 \begin{aligned}
 \|(x,v)-S_N(x,v)\|_{X\oplus_1V}
 &\le \|(x,v)-P_{j-1}(x,v)\|_{X\oplus_1V}\\
 &\quad+\|B_{j,k}\|\bigl(\|d_j\|_X+\|v_j\|_{V_j}\bigr).
 \end{aligned}
\]
As \(N\to\infty\), we have \(j\to\infty\). The first term tends to zero
by the convergence of \(P_{j-1}(x,v)\), and the second tends to zero since
\(\|d_j\|_X+\|v_j\|_{V_j}\to0\). Thus \(S_N(x,v)\to(x,v)\).
Applying \(P_j-P_{j-1}\) to any convergent basis expansion of \((x,v)\)
recovers the expansion of \((d_j,v_j)\) in the \(j\)th block, where the
coefficients are unique. This gives a Schauder basis of \(X\oplus_1V\).
Transporting it to \(X\) by the isomorphism above and then applying the Bohr
transform gives a Schauder basis of \(\mathcal H^1\).

\subsection{An explicit Schauder basis of \texorpdfstring{\(\mathcal H^1\)}{H\string^1}}\label{sec:explicit}

We construct the basis using the blocks
\((P_n-P_{n-1})X\), where the projections \(P_n\) are defined in
\eqref{eq:prime-block-projection}. The first three
vectors are \(1,2^{-s},25^{-s}\). Subsequent vectors are indexed by a block
number \(n\ge2\) and a matrix \(\beta\in\mathcal D_n\), with \(\mathcal D_n\) as in \eqref{eq:matrix-indices}. Formula~\eqref{eq:explicit-vector}
gives the vector and \eqref{eq:explicit-dual} gives its coefficient functional.
We use the binary order of Proposition~\ref{prop:prefix}.

\subsubsection{Finite-rank projections on \texorpdfstring{\(H^1(\mathbb T^\infty)\)}{the Hardy space on the infinite torus}}
Fix \(n\ge1\) and put \(d=2^n-1\) and \(K=2d\).
The auxiliary maps \(J,Q,A,B\) below depend on this choice of \(n\).
Group the first \(2d+1\) variables of \((t,z_1,z_2,\ldots)\), in order,
as \((\xi,a,v)\in\mathbb T^d\times\mathbb T\times\mathbb T^d\).
Let \(E_n\) be the span of \(\xi^\alpha\) with \(\alpha\in\mathbb N_0^d\) and \(|\alpha|\le K\).
Let \(A_d:X\to H^1(\mathbb T^d)\) be the section in the first \(d\) coordinates
defined by \eqref{eq:coordinate-section}.
Define \(J:E_n\to X\) and \(Q:X\to E_n\) by
\[
 \begin{aligned}
 &(Jh)(\xi,a,v)=a^K h(v_1/a,\ldots,v_d/a),\\
 &(Qf)(\xi)=\int_{\mathbb T^d}\int_{\mathbb T}\bar\lambda^K
   (A_{2d+1}f)(w,\lambda,\lambda\xi)\,dm(\lambda)\,dm(w).
 \end{aligned}
\]
Here \(\lambda\xi=(\lambda\xi_1,\ldots,\lambda\xi_d)\).
Integration in \(v\) for fixed \(a\) gives \(\|Jh\|_X=\|h\|_{L^1(\mathbb T^d)}\).
For \(Q\), integration of the bound for \(|(Qf)(\xi)|\) in \(\xi\) yields
\(\|Qf\|_{L^1(\mathbb T^d)}\le\|A_{2d+1}f\|_{L^1(\mathbb T^{2d+1})}\le\|f\|_X\).
On monomials, \(J\xi^\alpha=a^{K-|\alpha|}v^\alpha\) for \(|\alpha|\le K\).
The integral gives \(Q(a^{K-|\alpha|}v^\alpha)=\xi^\alpha\) for every
\(|\alpha|\le K\), and \(Q\) vanishes on every analytic monomial outside this family.
Hence \(QJ=I_{E_n}\), and \(A_dJ=0\)
because every monomial of \(Jh\) has degree \(K>0\) in \((a,v)\).

With \(b(x)=e^{-1/x}\) for \(x>0\) and \(b(x)=0\) for \(x\le0\), set
\[
 \theta(x)=\frac{\pi}{2}\frac{b(x-1)}{b(x-1)+b(2-x)} .
\]
Define \(A,B\) on analytic polynomials on \(\mathbb T^d\) by
\[
 A \xi^\alpha=\cos\theta(|\alpha|/d)\xi^\alpha,\qquad
 B \xi^\alpha=\sin\theta(|\alpha|/d)\xi^\alpha .
\]
Let \(\varphi\) be the even extension from \([0,\infty)\) of either \(\cos\theta\) or \(\sin\theta-1\). Then \(\varphi\in C_c^\infty(\mathbb R)\). With \(\check\varphi(u)=\int_{\mathbb R}\varphi(v)e^{2\pi iuv}\,dv\), set
\[
 \kappa_d(x)=\sum_{l\in\mathbb Z}d\check\varphi(d(x+l)),\qquad x\in\mathbb R/\mathbb Z.
\]
Integration over the translated unit intervals gives
\(\widehat{\kappa_d}(j)=\int_0^1\kappa_d(x)e^{-2\pi ijx}\,dx=\varphi(j/d)\) and
\(\|\kappa_d\|_{L^1(\mathbb R/\mathbb Z)}\le\|\check\varphi\|_{L^1(\mathbb R)}\).
For any analytic polynomial \(h\) on \(\mathbb T^d\),
\[
 \int_0^1\kappa_d(x)h(e^{-2\pi ix}\xi)\,dx
 =\sum_\alpha \varphi(|\alpha|/d)\widehat h(\alpha)\xi^\alpha .
\]
The two choices of \(\varphi\) give \(Ah\) and \((B-I)h\), respectively.
The integral has \(L^1(\mathbb T^d)\) norm at most
\(\|\check\varphi\|_{L^1(\mathbb R)}\|h\|_{L^1(\mathbb T^d)}\).
Thus \(A,B:H^1(\mathbb T^d)\to H^1(\mathbb T^d)\) are bounded
uniformly in \(n\). Since \(\cos\theta(x)=0\) for \(x\ge2\), \(A\) has
range in \(E_n\). Set
\begin{equation}\label{eq:prime-block-projection}
 \begin{aligned}
 &V_n:E_n\to X, &&V_n h=Ah+JBh,\\
 &L_n:X\to E_n, &&L_nf=A(A_df)+BQf,
 \end{aligned}
 \qquad P_n=V_nL_n:X\to X,\quad P_0=0 .
\end{equation}
Since \(A_dJ=0,\ QJ=I\), and \(Q\) vanishes on polynomials in \(\xi\),
\[
 L_nV_n h=A^2h+B^2h=h\qquad(h\in E_n).
\]
Thus \(P_n^2=P_n\) and \(\sup_n\|P_n\|\le\sup_n\|V_n\|\|L_n\|<\infty\).
For \(\alpha\in\mathbb N_0^d\) with \(|\alpha|\le d\), \(\theta(|\alpha|/d)=0\), so
\[
 P_n \xi^\alpha=\xi^\alpha,\qquad
 \widehat{P_n f}(\alpha)=\widehat f(\alpha)\quad(f\in X),
\]
where \(\alpha\) is extended by zeros. By the formulas for \(V_n,L_n\),
\(P_nf\) is a linear combination of monomials of total degree at most \(2d\)
in the first \(2d+1\) variables and depends only on the coefficients of \(f\)
at those monomials.
For \(l>n\) we have \(2^l-1\ge2d+1\), so \(P_l\) fixes these monomials
and preserves their coefficients. It follows that \(P_lP_n=P_nP_l=P_n\).
Set \(\Delta_n=P_n-P_{n-1}:X\to X\). Since \(P_n\) fixes every analytic
polynomial for all sufficiently large \(n\), density and uniform boundedness give
\(\sum_{j=1}^n\Delta_j f=P_nf\to f\) in \(X\) for every \(f\in X\).
We use the spaces \(\Delta_n X\) as the blocks of the explicit construction.
In the first block \(d=1\) and \(K=2\), so
\begin{equation}\label{eq:first-block-projection}
 (P_1f)(t,z)=\widehat f(0)+\widehat f(\alpha(2))\,t+\widehat f(\alpha(25))\,z_2^2 .
\end{equation}

\subsubsection{The finite matrices on each block}

For \(n\ge2\), take \(I_n:E_n\to W_n\) and \(M_n:W_n\to E_n\) as in
\eqref{eq:matrix-embedding} and \eqref{eq:matrix-left-inverse}, respectively, with the parameters in \eqref{eq:block-parameters}. Their norms are uniformly bounded in \(n\) by Lemma~\ref{lem:row-bound}. Set

\begin{equation}\label{eq:explicit-block-maps}
 \begin{aligned}
 &i_n:\Delta_n X\to W_n, &&i_nx=I_n L_nx,\\
 &r_n:W_n\to\Delta_n X, &&r_nw=(I-P_{n-1})V_n M_nw .
 \end{aligned}
\end{equation}
Since \(P_n V_n=V_n\), \(r_nw=\Delta_n V_nM_nw\).
For \(x\in\Delta_n X\), \(M_nI_n=I_{E_n}\) and \eqref{eq:prime-block-projection} give
\[
 r_ni_nx=(I-P_{n-1})V_nL_nx=(I-P_{n-1})P_nx=x .
\]
Put \(K_n=\ker r_n\) and \(T_n=I_{W_n}-i_nr_n\).
Since \(r_nT_n=0\) and \(T_n k=k\) for \(k\in K_n\), \(T_n\) is a
projection onto \(K_n\), with \(\|T_n\|\le1+\|i_n\|\|r_n\|\) uniformly in \(n\).
In particular \(w\mapsto(r_nw,T_nw)\) and
\((x,k)\mapsto i_nx+k\) are mutually inverse between \(W_n\) and
\(\Delta_n X\oplus_1 K_n\).

\subsubsection{Returning to the original variables}
Fix \(n\ge2\), set \(d=2^n-1\), and take \(r,a_0\) as in \eqref{eq:block-parameters}.
For \(z=(z_k)_{k\ge1}\in\mathbb T^\infty\), define the matrices
\(u=(u_{ij})\), \(v=(v_{ij})\), the vector \(z'=(z'_1,\ldots,z'_d)\), and \(\ell\) by
\[
 \begin{aligned}
 &u_{ij}=z_{(i-1)r+j}, &&v_{ij}=z_{dr+(i-1)r+j}
       &&(1\le i\le d,\ 1\le j\le r),\\
 &z'_i=z_{2dr+i}, &&\ell=z_{2dr+d+1} &&(1\le i\le d).
 \end{aligned}
\]
Define \(G_n:W_n\to H^1(\mathbb T^\infty)\) by
\begin{equation}\label{eq:prime-monomial-embedding}
 G_n w^\beta=\left(\prod_{i=1}^d\prod_{j=1}^r u_{ij}^{\beta_{ij}}v_{ij}^{1-\beta_{ij}}\right)(z')^{\nu(\beta)}
                  \ell^{K-|\nu(\beta)|}.
\end{equation}
Substituting \(\nu_i(\beta)=\sum_j\beta_{ij}-a_0\) gives, for \(F\in W_n\),
\[
 (G_nF)(z)=\left(\prod_{i,j}v_{ij}\right)
 \left(\prod_i (z'_i)^{-a_0}\right)\ell^{K+da_0}
 F\!\left(\left(\frac{u_{ij}z'_i}{v_{ij}\ell}\right)_{i,j}\right).
\]
Integrating first in \(u\) gives \(\|G_nF\|_{H^1(\mathbb T^\infty)}=\|F\|_{W_n}\).
Write \(G_nw^\beta=z^{\gamma(n,\beta)}\), where \(\gamma_j(n,\beta)\)
is the exponent of \(z_j\) in \eqref{eq:prime-monomial-embedding}. Define
\[
 G_n^\flat:H^1(\mathbb T^\infty)\to W_n,\qquad
 G_n^\flat f=\sum_{\beta\in\mathcal D_n}\widehat f(\gamma(n,\beta))w^\beta .
\]
By \eqref{eq:prime-monomial-embedding}, the exponent of \(u_{ij}=z_{(i-1)r+j}\) in \(G_nw^\beta\) is
\(\gamma_{(i-1)r+j}(n,\beta)=\beta_{ij}\). Thus distinct matrices give distinct monomials, and for \(\beta,\delta\in\mathcal D_n\),
\[
 \widehat{G_nw^\beta}\bigl(\gamma(n,\delta)\bigr)=1_{\{\delta=\beta\}},
 \qquad G_n^\flat G_nw^\beta=w^\beta .
\]
These monomials \(w^\beta\) span \(W_n\), so \(G_n^\flat G_n=I_{W_n}\).
To see that \(\|G_n^\flat\|\le1\), first average over the variables \(z_k\) with \(k>2dr+d+1\).
In the remaining coordinates apply \eqref{eq:affine-fourier-projection} with the affine map
\[
 \beta\longmapsto\left((1-\beta_{ij})_{i,j},\,
              (\textstyle\sum_j\beta_{ij}-a_0)_i,\,
              K+da_0-\textstyle\sum_{i,j}\beta_{ij}\right).
\]
For \(\beta\in\mathbb N_0^{d\times r}\), all entries of this tuple are nonnegative exactly when
\[
 \beta_{ij}\le1\quad(1\le i\le d,\ 1\le j\le r),\qquad
 \sum_{j=1}^r\beta_{ij}\ge a_0\quad(1\le i\le d),\qquad
 \sum_{i=1}^d\sum_{j=1}^r\beta_{ij}\le K+da_0 .
\]
By \eqref{eq:matrix-indices}, these conditions mean \(\beta\in\mathcal D_n\).
The projection thus retains exactly the monomials \(z^{\gamma(n,\beta)}\), \(\beta\in\mathcal D_n\), and equals \(G_nG_n^\flat\). Since \(G_n\) is isometric,
\(\|G_n^\flat f\|_{W_n}=\|G_nG_n^\flat f\|_{H^1(\mathbb T^\infty)}\le\|f\|_{H^1(\mathbb T^\infty)}\).
Put \(q_n=T_nG_n^\flat:H^1(\mathbb T^\infty)\to K_n\). Then \(q_nG_n k=k\) for \(k\in K_n\).

Order the pairs \((b,n)\), \(b\ge0,\ n\ge2\), by increasing \(b+n\),
breaking ties by increasing \(b\). For the \(j\)th pair, starting at \(j=1\), set
\(m_j=2^{2j+6}\), \(\zeta_j=\exp(i2^{-j})\), and
\begin{equation}\label{eq:fejer-polynomials}
 h_j(t)=\sum_{k=1-m_j}^{m_j-1}
 (1-|k|/m_j)\bar\zeta_j^{\,k}t^{4m_j+k},\qquad
 k_j=\frac{3m_j}{2m_j^2+1}h_j .
\end{equation}
The Fej\'er kernel formula \cite[Chapter~2, Lemma~5.1]{SteinShakarchi} gives
\[
 |h_j(t)|=m_j^{-1}\left|\sum_{l=0}^{m_j-1}(\bar\zeta_jt)^l\right|^2,
 \qquad \|h_j\|_{L^1(\mathbb T)}=1 .
\]
By \eqref{eq:fejer-polynomials}, \(\widehat h_j(k)=\widehat k_j(k)=0\) for \(k\notin[3m_j+1,5m_j-1]\). These intervals are disjoint since \(m_{j+1}=4m_j\).
Parseval's identity gives
\begin{equation}\label{eq:fejer-biorthogonality}
 \|h_j\|_{L^2(\mathbb T)}^2=1+2\sum_{k=1}^{m_j-1}(1-k/m_j)^2
       =\frac{2m_j^2+1}{3m_j},\qquad
 \int_{\mathbb T}h_j\overline{k_i}\,dm=1_{\{i=j\}}.
\end{equation}
Writing \(t=\zeta_j e^{i\theta}\), \(-\pi\le\theta\le\pi\), the geometric sum satisfies
\(\big|\sum_{l=0}^{m_j-1}e^{il\theta}\big|
 \le\min\{m_j,2/|1-e^{i\theta}|\}\) for \(\theta\ne0\). Hence
\[
 |k_j(\zeta_j e^{i\theta})|\lesssim(1+m_j|\theta|)^{-2} .
\]
The arcs of angular radius \(2^{-j-3}\) about \(\zeta_j\) are disjoint.
Outside the \(j\)th arc, \(m_j|\theta|\ge2^{j+3}\).
Consequently \(\sup_t\sum_j|k_j(t)|\lesssim 1+\sum_{j\ge1}2^{-2j}<\infty\).

Let \(Z=(\bigoplus_{b\ge0,n\ge2}K_n)_1\), with \(\|v\|_Z=\sum_{b,n}\|v_{b,n}\|_{W_n}\).
If \((b,n)\) is the \(j\)th pair in the ordering above, we write \(h_{b,n}=h_j\) and \(k_{b,n}=k_j\)
for the polynomials in \eqref{eq:fejer-polynomials}, so that their subscripts match the summand of \(Z\) indexed by \((b,n)\).
Define \(S:Z\to X\)
and \(R:X\to Z\) by
\begin{equation}\label{eq:explicit-sum-maps}
 (Sv)(t,z)=\sum_{b,n} h_{b,n}(t)\,(G_n v_{b,n})(z),\qquad
 (Rf)_{b,n}=q_n\!\left(\int_{\mathbb T} f(t,\cdot)\overline{k_{b,n}(t)}\,dm(t)\right).
\end{equation}
Since \(G_n\) is isometric and \(\|h_j\|_1=1\), the series for \(Sv\)
converges absolutely with \(\|Sv\|_1\le\|v\|_Z\). Fubini's theorem gives
\[
 \|Rf\|_Z\le \sup_n\|q_n\|\iint |f(t,z)|\sum_j|k_j(t)|\,dm(t)\,dm(z)
 \le \left(\sup_n\|q_n\|\right)\left(\sup_t\sum_j|k_j(t)|\right)\|f\|_1 .
\]
By \eqref{eq:fejer-biorthogonality} and \(q_nG_n|_{K_n}=I_{K_n}\), we have \(RS=I_Z\).
For \(v\in Z\) define \(\sigma,\sigma_-:Z\to Z\) by
\[
 (\sigma v)_{0,n}=0,\qquad(\sigma v)_{b+1,n}=v_{b,n},\qquad
 (\sigma_-v)_{b,n}=v_{b+1,n}.
\]
Define \(\iota_0:(\bigoplus_{n\ge2}K_n)_1\to Z\) by
\((\iota_0u)_{b,n}=1_{\{b=0\}}u_n\).
With \(\mathcal T=I-SR+S\sigma R:X\to X\) we obtain inverse isomorphisms
\begin{equation}\label{eq:explicit-isomorphism}
\begin{aligned}
 &X\oplus_1(\textstyle\bigoplus_{n\ge2}K_n)_1\longrightarrow X,
 &&(x,v)\longmapsto\mathcal T x+S\iota_0 v,\\
 &X\longrightarrow X\oplus_1(\textstyle\bigoplus_{n\ge2}K_n)_1,
 &&f\longmapsto\big((I-SR+S\sigma_-R)f,((Rf)_{0,n})_{n\ge2}\big).
\end{aligned}
\end{equation}
Indeed each \(f\in X\) has the unique form \(f=(I-SR)f+S(Rf)\)
with \((I-SR)f\in\ker R\). The first map fixes the \(\ker R\)
component and sends \((Rx,v)\) to \(\sigma Rx+\iota_0v\) in \(Z\).
Thus \((\sigma Rx+\iota_0v)_{0,n}=v_n\) for \(n\ge2\), and \(\sigma_-(\sigma Rx+\iota_0 v)=Rx\).

\subsubsection{Individual basis vectors and their coefficients}
For \(n\ge2\), \(\alpha\in\mathbb N_0^d\) and \(|\alpha|\le K\), put
\(Y_{n,\alpha}=\Delta_n V_n\xi^\alpha\) and
\[
 C_{n,\alpha}=\mathcal B
  \big(\mathcal T Y_{n,\alpha}-h_{0,n}G_n i_nY_{n,\alpha}\big).
\]
In the variables \((t,z_1,z_2,\ldots)\), the Bohr transform \eqref{eq:bohr-transform} gives
\[
 \mathcal B(t^k z^\gamma)(s)
 =\bigg(2^k\prod_{j\ge1}p_{j+1}^{\gamma_j}\bigg)^{-s}.
\]
For \(\beta\in\mathcal D_n\), put \(m_{n,\beta}=\prod_{j\ge1}p_{j+1}^{\gamma_j(n,\beta)}\).
With \(t_\beta\) as in \eqref{eq:matrix-left-inverse}, the Dirichlet basis vector indexed by \((n,\beta)\) is
\begin{equation}\label{eq:explicit-vector}
 B_{n,\beta}(s)=m_{n,\beta}^{-s}h_{0,n}(2^{-s})
                   +t_\beta C_{n,\nu(\beta)}(s).
\end{equation}
Indeed, the \(\Delta_n X\) and \(K_n\) components of \(w^\beta\)
in \eqref{eq:explicit-block-maps} are
\[
 x=r_nw^\beta=t_\beta Y_{n,\nu(\beta)},\qquad
 k=T_nw^\beta=w^\beta-i_nx .
\]
Placing \(k\) in the \(n\)th summand of \((\bigoplus_{l\ge2} K_l)_1\)
gives the image \(\mathcal T x+h_{0,n}G_n(w^\beta-i_nx)\) under \eqref{eq:explicit-isomorphism}.
Applying the Bohr transform \eqref{eq:bohr-transform} yields \eqref{eq:explicit-vector}.

For a monomial \(t^k z^\gamma\), at most one of the intervals \([3m_j+1,5m_j-1]\) contains \(k\).
If \(k\) belongs to the interval for the pair \((b,l)\), form the matrix
\(\delta_{ic}=\gamma_{(i-1)r+c}\), \(1\le i\le d\), \(1\le c\le r\), with \(d,r\) as in block \(l\).
If \(\delta\in\mathcal D_l\) and \(\gamma=\gamma(l,\delta)\), then \eqref{eq:explicit-sum-maps} gives
\[
 \mathcal T(t^kz^\gamma)=t^kz^\gamma+
    \overline{\widehat{k_{b,l}}(k)}\,(h_{b+1,l}-h_{b,l})G_l T_l w^\delta .
\]
In all other cases \(\mathcal T(t^kz^\gamma)=t^kz^\gamma\).
Thus \eqref{eq:explicit-vector} requires only finite sums.

As in \eqref{eq:binary-order}, order \(\mathcal D_n\)
by increasing \(\sum_{i,j}2^{(i-1)r+j-1}\beta_{ij}\).
By \eqref{eq:matrix-count}, the cumulative block lengths are determined by
\[
 N_1=3,\qquad
 N_n=N_{n-1}+|\mathcal D_n|,\qquad
 |\mathcal D_n|=\sum_{\substack{\alpha\in\mathbb N_0^d\\|\alpha|\le K}}
                       \prod_{i=1}^d\binom r{a_0+\alpha_i}.
\]
For \(k>3\), choose \(n\) so that \(N_{n-1}<k\le N_n\), and let \(\beta\)
be the matrix at position \(k-N_{n-1}\) in the above ordering of \(\mathcal D_n\).
Then \(B_{n,\beta}\) is the \(k\)th basis vector.
In block \(2\), the least matrix \(\beta_*\) has ones in the
first \(r/2\) columns of each row and zeros elsewhere. Hence
\(\nu(\beta_*)=0\), and \(Y_{2,0}=(I-P_1)V_2 1=0\).
Substituting in \eqref{eq:explicit-vector} and using \eqref{eq:fejer-polynomials} gives the fourth basis vector
\[
 B_{2,\beta_*}(s)=\sum_{k=-255}^{255}(1-|k|/256)e^{-ik/2}
                  (2^{1024+k}m_{2,\beta_*})^{-s}.
\]

For \(F\in\mathcal H^1\) put \(f=\mathcal B^{-1} F\) and
\(g=(I-SR+S\sigma_-R)f\). The \(n\)th block of the inverse image in
\eqref{eq:explicit-isomorphism} is \((\Delta_n g,(Rf)_{0,n})\) for \(n\ge2\).
Under the inverse block map this pair corresponds to \(i_n\Delta_n g+(Rf)_{0,n}\in W_n\).
Its Fourier coefficient at \(\beta\) gives
\begin{equation}\label{eq:explicit-dual}
 c_{n,\beta}(F)=\int_{\mathbb T^{dr}}
 \big(i_n\Delta_n g+(Rf)_{0,n}\big)(w)\,w^{-\beta}\,dm(w).
\end{equation}
The monomials \(1,t,z_2^2\) span \(P_1 X\) by \eqref{eq:first-block-projection}.
Since \(\widehat{k_{b,n}}(0)=\widehat{k_{b,n}}(1)=0\) for every \((b,n)\),
\eqref{eq:explicit-sum-maps} gives \(R1=Rt=R(z_2^2)=0\), so \(\mathcal T\) fixes these monomials.
Likewise \(\widehat{h_{b,n}}(0)=\widehat{h_{b,n}}(1)=0\) gives \(P_1S=0\), hence \(P_1g=P_1f\). Thus for
\(F(s)=\sum_{m\ge1}a_m m^{-s}\) the first three basis vectors are
\(1,2^{-s},25^{-s}\) and their respective coefficients are \(a_1,a_2,a_{25}\).

To compute the remaining coefficients, write \(\widehat f(k,\gamma)\)
for the coefficient of \(t^kz^\gamma\), so that
\(\widehat f(k,\gamma)=a_{2^k\prod_{j\ge1}p_{j+1}^{\gamma_j}}\).
If \((b,n)\) is the \(j\)th pair, then
\[
 (Rf)_{b,n}=
 T_n\left(\sum_{\delta\in\mathcal D_n}
       \sum_{k=3m_j+1}^{5m_j-1}
       \overline{\widehat{k_{b,n}}(k)}\,\widehat f(k,\gamma(n,\delta))\,w^\delta\right).
\]
Thus each \((Rf)_{b,n}\) uses finitely many coefficients of \(F\).
If \(3m_j+1\le k\le 5m_j-1\) and \((b,l)\) is the \(j\)th pair, then
\(g-f=S(\sigma_--I)Rf\) gives
\[
 \widehat{g-f}(k,\gamma)
 =\widehat{h_{b,l}}(k)\,
   \widehat{G_l\big((Rf)_{b+1,l}-(Rf)_{b,l}\big)}(\gamma).
\]
If \(k\) lies outside all these intervals, then \(\widehat{g-f}(k,\gamma)=0\).
Since \(\Delta_n g\) uses only finitely many coefficients of \(g\), this
also proves the finiteness assertion for \eqref{eq:explicit-dual}.
If \(F\) is a Dirichlet polynomial, \(f=\mathcal B^{-1}F\) is a polynomial in \(t,z\).
Choose \(M\) so that \(\widehat f(k,\gamma)=0\) for every \(k>M\) and every \(\gamma\).
Then \((Rf)_{b,n}=0\) whenever \(3m_j+1>M\), where \(j\) indexes \((b,n)\).
Hence \(Rf\) is finitely supported and \(g\) is a polynomial. For all sufficiently
large \(n\), both \(\Delta_n g\) and \((Rf)_{0,n}\) vanish. The inverse
block maps therefore express \(F\) as a finite linear combination of the
stated basis vectors.

In \(X\oplus_1(\bigoplus_{l\ge2}K_l)_1\), the partial sum through block \(n\) maps
\((x,(v_l)_{l\ge2})\) to \((P_nx,(1_{\{l\le n\}}v_l)_{l\ge2})\). These operators have
norm at most \(\max\{1,\sup_n\|P_n\|\}\).
Within block \(n\ge2\), the map \(w\mapsto(r_nw,T_nw)\) and its inverse
have uniformly bounded norms, so \eqref{eq:binary-prefix-bound} gives a
uniform bound on the partial sums in \(\Delta_n X\oplus_1 K_n\).
Formula~\eqref{eq:concatenated-partial-sum} therefore bounds every partial
sum in the stated order. The isomorphism \eqref{eq:explicit-isomorphism}
transfers this bound to \(X\). Since these partial sums eventually fix each
polynomial, they converge in norm on \(X\) by density. Uniqueness follows
by applying the inverse in \eqref{eq:explicit-isomorphism}, then the block
projections and the coefficient functionals on \(W_n\).

This completes the proof of Theorem~\ref{thm:basis}.\qed

\medskip
\noindent\textbf{AI Assistance Disclosure.}
OpenAI's GPT-6 Astra provided mathematical suggestions that contributed to the development of parts of the proofs of Theorem~\ref{thm:retract} and Lemma~\ref{lem:row-bound}.
The explicit construction and representation of the Schauder basis presented in Section~\ref{sec:explicit} were also developed with the assistance of GPT-6 Astra.
The authors take full responsibility for the mathematical correctness, validity, and integrity of all results and proofs in this paper.

\begingroup
\footnotesize
\newcommand{\authorcontact}[1]{%
  \par\medskip
  \textsc{School of Mathematical Sciences, Dalian University of Technology,
  Dalian, Liaoning 116024, P. R. China}\\*
  \textit{Email address:} \href{mailto:#1}{\nolinkurl{#1}}\par}
\authorcontact{503304283@qq.com}
\authorcontact{lyfdlut@dlut.edu.cn}
\authorcontact{zuchao@dlut.edu.cn}
\endgroup
\end{document}